\documentclass[a4paper]{amsart}
\usepackage{graphics}
\usepackage{amsmath}
\usepackage{latexsym}
\usepackage{amssymb}
\usepackage[all]{xy}
\xyoption{matrix}
\xyoption{arrow}

\begin{document}

\newcommand{\lcm}{\operatorname{lcm}\nolimits}
\newcommand{\coker}{\operatorname{coker}\nolimits}
\renewcommand{\th}{\operatorname{th}\nolimits}
\newcommand{\rej}{\operatorname{rej}\nolimits}
\newcommand{\extto}{\xrightarrow}
\renewcommand{\mod}{\operatorname{mod}\nolimits}
\newcommand{\Sub}{\operatorname{Sub}\nolimits}
\newcommand{\ind}{\operatorname{ind}\nolimits}
\newcommand{\Fac}{\operatorname{Fac}\nolimits}
\newcommand{\add}{\operatorname{add}\nolimits}
\newcommand{\Hom}{\operatorname{Hom}\nolimits}
\newcommand{\Rad}{\operatorname{Rad}\nolimits}
\newcommand{\RHom}{\operatorname{RHom}\nolimits}
\newcommand{\uHom}{\operatorname{\underline{Hom}}\nolimits}
\newcommand{\End}{\operatorname{End}\nolimits}
\renewcommand{\Im}{\operatorname{Im}\nolimits}
\newcommand{\Ker}{\operatorname{Ker}\nolimits}
\newcommand{\Coker}{\operatorname{Coker}\nolimits}
\newcommand{\Ext}{\operatorname{Ext}\nolimits}
\newcommand{\op}{{\operatorname{op}}}
\newcommand{\Ab}{\operatorname{Ab}\nolimits}
\newcommand{\id}{\operatorname{id}\nolimits}
\newcommand{\pd}{\operatorname{pd}\nolimits}
\newcommand{\ql}{\operatorname{q.l.}\nolimits}
\newcommand{\rank}{\operatorname{rank}\nolimits}
\newcommand{\A}{\operatorname{\mathcal A}\nolimits}
\newcommand{\C}{\operatorname{\mathcal C}\nolimits}
\newcommand{\D}{\operatorname{\mathcal D}\nolimits}
\newcommand{\E}{\operatorname{\mathcal E}\nolimits}
\newcommand{\X}{\operatorname{\mathcal X}\nolimits}
\newcommand{\Y}{\operatorname{\mathcal Y}\nolimits}
\newcommand{\F}{\operatorname{\mathcal F}\nolimits}
\newcommand{\Z}{\operatorname{\mathbb Z}\nolimits}
\renewcommand{\P}{\operatorname{\mathcal P}\nolimits}
\newcommand{\T}{\operatorname{\mathcal T}\nolimits}
\newcommand{\G}{\Gamma}
\renewcommand{\L}{\Lambda}
\newcommand{\bdot}{\scriptscriptstyle\bullet}
\renewcommand{\r}{\operatorname{\underline{r}}\nolimits}
\newtheorem{lem}{Lemma}[section]
\newtheorem{prop}[lem]{Proposition}
\newtheorem{cor}[lem]{Corollary}
\newtheorem{thm}[lem]{Theorem}
\newtheorem*{thmA}{Theorem}
\newtheorem*{thmB}{Theorem}
\newtheorem{rem}[lem]{Remark}
\newtheorem{defin}[lem]{Definition}
\newtheorem{example}[lem]{Example}


\title[Denominators of cluster variables]{Denominators of cluster
variables}

\author[Buan]{Aslak Bakke Buan}
\address{Institutt for matematiske fag\\
Norges teknisk-naturvitenskapelige universitet\\
N-7491 Trondheim\\
Norway}
\email{aslakb@math.ntnu.no}

\author[Marsh]{Robert J. Marsh}
\address{Department of Pure Mathematics \\
University of Leeds \\
Leeds \\
LS2 9JT \\
England
}
\email{marsh@maths.leeds.ac.uk}

\author[Reiten]{Idun Reiten}
\address{Institutt for matematiske fag\\
Norges teknisk-naturvitenskapelige universitet\\
N-7491 Trondheim\\
Norway}
\email{idunr@math.ntnu.no}


\begin{abstract}
Associated to any acyclic cluster algebra is a corresponding triangulated
category known as the cluster category. It is known that there is a
one-to-one correspondence between cluster variables in the cluster algebra
and exceptional indecomposable objects in the cluster category inducing
a correspondence between clusters and cluster-tilting objects.

Fix a cluster-tilting object $T$ and a corresponding initial cluster.
By the Laurent phenomenon, every cluster variable can be written as a
Laurent polynomial in the initial cluster.
We give conditions on $T$ equivalent to the fact that
the denominator in the reduced form for every cluster variable in the
cluster algebra has exponents given by the dimension vector of the
corresponding module over the endomorphism algebra of $T$.
\end{abstract}

\thanks{The authors were supported by Storforsk grant no.\ 167130 from
the Norwegian Research Council.}

\maketitle

\section*{Introduction}

Cluster algebras were introduced by Fomin-Zelevinsky
\cite{fominzelevinsky02}, and have been influential in many
settings including algebraic combinatorics, Lie theory, Poisson
geometry, Teichm{\"u}ller theory, total positivity and quiver
representations.

The first link to quiver representations and tilting theory for
algebras was given in \cite{mrz03}. Cluster categories were
defined in \cite{bmrrt06}, giving a categorical model for cluster
combinatorics. For type $A$, an independent approach was given in
\cite{ccs06}.

A cluster algebra is in essence a commutative ring with a
distinguished countable family of generators, called \emph{cluster
variables}. The cluster variables can be grouped into overlapping
sets of equal finite size $n$, called \emph{clusters}. Any given
cluster variable can be uniquely expressed as a rational function
in the elements of a fixed initial cluster. By the Laurent
phenomenon, \cite[3.1]{fominzelevinsky02}, these rational
functions are actually Laurent polynomials with integer
coefficients.

Corresponding to each cluster there is a quiver $Q$, with $n$
vertices. The cluster algebra is said to be acyclic if it admits
a cluster for which the corresponding quiver $Q$ has no oriented cycles.
In this case there is a corresponding finite dimensional hereditary
algebra, the path algebra of $kQ$ (for some field $k$) giving rise to a cluster
category.

In \cite{bmrt} a surjective map $\alpha$ from the cluster
variables of an acyclic cluster algebra to the exceptional indecomposable
objects in the corresponding cluster category was given.
In~\cite{calderokellerB}
a map was defined in the opposite direction, and it was shown that
this map is a bijection with $\alpha$ as inverse (see also~\cite{bckmrt}).
Our exposition is made in such a way that it is easy to see where the fact
that $\alpha$ is a bijection is used.

An interesting property of the map $\alpha$ from \cite{bmrt} is
the following. Choose a cluster whose corresponding quiver is acyclic.
If a cluster variable (not from this cluster) is expressed as a Laurent
polynomial in the cluster variables from this cluster there is a natural
interpretation of the monomial in the denominator in terms of the
composition factors of its image under $\alpha$ (regarded as a
$kQ$-module). This generalized earlier results in this
direction~\cite{fominzelevinsky03,ccs,calderokellerA,reitentodorov}.

It is interesting to ask whether for an arbitrary choice of initial cluster it 
is possible to interpret the monomial in the denominator of a
cluster variable in terms of the composition factors of the module
(over an appropriate cluster-tilted algebra) corresponding to the
image of the cluster variable under $\alpha$. The cluster-tilted
algebra involved should be the endomorphism algebra of the
cluster-tilting object $\tau T$ obtained by forming the direct sum of
the images of the fixed initial cluster under $\alpha$.
We remark that such an interpretation has been found in the Dynkin case
in~\cite{ccs,reitentodorov}.

In this paper we give a precise answer to this question for
the general case by giving
necessary and sufficient conditions on the cluster-tilting object $T$ for
this to hold. For a given acyclic cluster algebra, these
conditions hold for all cluster-tilting objects $T$ in the corresponding cluster category 
if and only if the algebra is of
finite cluster type or of rank 2. In the tame case we show that a
cluster-tilting object $T$ satisfies the conditions if and only if no
indecomposable regular summand of $T$ has quasilength exactly one less than the rank of
the tube in which it lies, and that this is equivalent to $\End_{\C}(T_i)\simeq k$ for
each indecomposable summand of $T$. In fact this last condition is
necessary for any path algebra $kQ$, and we conjecture that it is also sufficient in general.

We would like to thank Otto Kerner for answering several questions on wild algebras, and especially
for assisting us with the proof of Proposition \ref{p:endnotk}.

\section{Background and Main results} \label{s:mainresults}

Let $Q$ be a finite connected acyclic quiver and $k$ an algebraically closed field. Let $H=kQ$
be the corresponding (finite dimensional) path algebra and
$D^b(kQ)$ the bounded derived category of finite dimensional
$kQ$-modules. It has autoequivalences $\tau$ (the Auslander-Reiten
translate) and $[1]$ (the shift). Let
$\mathcal{C}=D^b(kQ)/\tau^{-1}[1]$ be the corresponding cluster
category. We shall regard modules over $kQ$ as objects of $\C$
(via the natural embedding of the module category over $kQ$ in $D^b(kQ)$).
For $i$ a vertex of $Q$ let $P_i$ denote the corresponding
indecomposable projective $kQ$-module.

Let $\mathcal{A}=\mathcal{A}(Q)\subseteq
\mathbb{F} =\mathbb{Q}(x_1,x_2,\ldots ,x_n)$ be the (acyclic,
coefficient-free) cluster algebra defined using the initial seed
$(\mathbf{x},Q)$, where $\mathbf{x}$ is the transcendence basis
$\{x_1,x_2,\ldots ,x_n\}$ of $\mathbb{F}$.

For an object $M$ of $\mathcal{C}$, let
$c_M=\prod_{i=1}^n x_i^{\dim \Hom_{\C}(P_i,M)}$.
For a polynomial $f=f(z_1,z_2,\ldots ,z_n)$, we say that $f$ satisfies
the \emph{positivity condition} if $f(e_i)>0$ where $e_i=(1,\ldots ,1,0,1,
\ldots ,1)$ (with a $0$ in the $i$th position) for $i=1,2,\ldots ,n$.
This condition was crucial in the investigations in \cite{bmrt},
and also plays an essential role here.

We first recall the following results. 

\begin{thm}\cite[2.2]{bmrt} \label{bmrtfirst}
Every cluster variable of $\mathcal{A}$ is either of the form $f/c_M$
(for some exceptional $kQ$-module $M$) or of the form
$fx_i$ for some $i$, where
in either case $f=f(x_1,x_2,\ldots ,x_n)$ satisfies the positivity
condition. Every indecomposable exceptional $kQ$-module arises in this way, and if
$M$ and $M'$ are two such modules satisfying $c_M=c_{M'}$ then $M\simeq M'$.
\end{thm}

A map $\alpha$ is then defined from cluster variables of $\mathcal{A}$ to
(isomorphism classes of) exceptional indecomposable objects of $\C$ in the
following way. If $f/c_M$ is a cluster variable for some exceptional
indecomposable $kQ$-module $M$ then $\alpha(f/c_M)=M$. If $fx_i$ is a
cluster variable then $\alpha(fx_i)=P_i[1]$.

Recall that a \emph{tilting seed} for $\mathcal{A}$ is a pair $(T,Q')$ where
$T$ is a cluster-tilting object of $\C$ and $Q'$ is the quiver of
$\End_{\C}(T)^{\op}$.

\begin{thm}\cite[2.3]{bmrt} \label{bmrtsecond}
The map $$\alpha \colon \{\text{cluster variables of $\mathcal{A}(Q)$} \} \to
  \{\text{indecomposable exceptional objects in }\C\}$$ is
  surjective. It induces a surjective map $\bar{\alpha} \colon \{
  \text{clusters}\} \to \{ \text{cluster-tilting objects}\}$, and
a surjective map $\widetilde{\alpha} \colon \{
  \text{seeds}\} \to \{ \text{tilting seeds}\}$, preserving quivers.
\end{thm}

We also have the following result of \cite{calderokellerA}, see also \cite{bckmrt}.

\begin{thm} \label{calderokellerresult}
There is a bijection $\beta:M\rightarrow X_M$ from indecomposable exceptional
objects of $\C$ to cluster variables of $\mathcal{A}$ such that for
any indecomposable exceptional $kQ$-module $M$, we have $X_M=f/c_M$
as an irreducible fraction of integral polynomials in the $x_i$.
\end{thm}

In particular $\alpha$ is an inverse of $\beta$ and for the cluster
variables $fx_i$ considered in Theorem~\ref{bmrtfirst}, we can only have
$f=1$. 

Let $\Gamma$ be a quiver mutation-equivalent to $Q$. Then there is
a seed $(\mathbf{y},\Gamma)$ of $\mathcal{A}$, where
$\mathbf{y}=\{y_1,y_2,\ldots ,y_n\}$ is a transcendence basis of
$\mathbb{F}$. For $i=1,2,\ldots ,n$ let
$T_i=\tau^{-1}\alpha(y_i)$ so that we have $\alpha(y_i)=\tau
T_i$. It follows from Theorem~\ref{bmrtsecond} that $\amalg_{i=1}^n \tau
T_i$ is a cluster-tilting object in $\C$ and that the quiver of
its endomorphism algebra in $\C$ is $\Gamma$. Hence the same holds
for $T=\amalg_{i=1}^n T_i$.

We can ask whether results analogous to the above hold for cluster
variables when they are expressed in terms of the initial cluster
$\{y_1,y_2,\ldots ,y_n\}$. For $M$ an object in $\C$, let
$t_M=\prod_{i=1}^n y_i^{\dim \Hom_{\C}(T_i,M)}$.

An interesting open question is whether $t_M$ (or, equivalently,
the dimensions of the spaces $\Hom_{\C}(T_i,M)$ for $i=1,2,\ldots ,n$)
determines $M$ uniquely. This is true if $T=kQ$ is the direct sum of
the non-isomorphic indecomposable projective modules but is not known in general.
In particular, this means that an approach following~\cite{bmrt} exactly
in this more general context is not possible.

\begin{defin} \label{typestar}
Let $x$ be a cluster variable of $\mathcal{A}$. We say that $x$ 
expressed in terms of the cluster $\mathbf{y}$ has a {\em $T$-denominator}
if either: \\
(I) We have that $\alpha(x)=M$ for some exceptional indecomposable object
$M$ of $\C$ not isomorphic to $\tau T_i$ for any $i$, and $x=f/t_M$, or \\
(II) We have that $\alpha(x)=\tau T_i$ for some $i$ and $x=fy_i$. \\
In either case, $f$ is a polynomial in the $y_i$ satisfying the
positivity condition.
\end{defin}

Our first main result is the following:

\begin{thm} \label{main1}
Let $Q$ be a finite quiver with no oriented cycles, and let $k$
be an algebraically closed field. Let $\C$ be the cluster category associated to $kQ$, and
let $T=\amalg_{i=1}^n T_i$ be a cluster-tilting object in $\C$. Let
$\mathcal{A}=\mathcal{A}(Q)$ be the cluster algebra associated to $Q$.
Then: 
\begin{itemize}
\item[(a)] If no indecomposable direct summand of $T$ is regular then every
cluster variable of $\mathcal{A}$ has a $T$-denominator.
\item[(b)] If every cluster variable of $\mathcal{A}$ has a $T$-denominator, then $\End_{\C}(T_i)\simeq k$ for all $i$.
\end{itemize}
\end{thm}

We note that, as in Theorem~\ref{bmrtfirst},
the bijectivity of $\alpha$ can be used
to show that in case (II) of Definition \ref{typestar}, only $f = 1$ occurs.
The proof of Theorem~\ref{main1} will be completed at the end of
Section~\ref{s:thepreprojectivecase}.
In the special case where $kQ$ is a tame algebra, we have the following:

\begin{thm} \label{main2}
Suppose that we are in the situation of Theorem~\ref{main1} and that,
in addition, $kQ$ is a tame algebra.
Then the following are equivalent: 
\begin{itemize}
\item[(a)] Every cluster variable of $\mathcal{A}$ has a $T$-denominator. 
\item[(b)] No regular summand $T_i$ of quasi-length $r-1$ lies in a tube of
rank $r$. 
\item[(c)] For all $i$, $\End_{\C}(T_i)\simeq k$.
\end{itemize}
\end{thm}

The proof of Theorem~\ref{main2} will be completed at the end of
Section~\ref{s:thetamecase}.
As a consequence of our main results, we have in the general case:

\begin{cor} \label{main3}
Let $Q$ be a finite quiver with no oriented cycles. 
Let $\C$ be the cluster category associated to $kQ$.
Let $\mathcal{A}=\mathcal{A}(Q)$ be the cluster algebra associated to $Q$.
Then every cluster variable of $\mathcal{A}$ has a $T$-denominator for
every cluster-tilting object $T$ if and only if $Q$ is Dynkin or has
exactly two vertices.
\end{cor}

The proof of Corollary~\ref{main3} will be completed at the end of
Section~\ref{s:thepreprojectivecase}.

\section{Conditions on a cluster-tilting object} \label{s:conditions}

In this section we consider certain conditions on a cluster-tilting
object which are equivalent to all cluster variables in
$\mathcal{A}$ having a $T$-denominator.
Fix an almost complete (basic) cluster-tilting object $\overline{T}'$
in $\mathcal{C}$. Let $M,M^{\ast}$ be the two complements of $\overline{T}'$,
so that $T'=\overline{T}'\amalg M$ and $T''=\overline{T}'\amalg M^{\ast}$ are
cluster-tilting objects (see~\cite[5.1]{bmrrt06}).
Let
\begin{equation}
M^{\ast} \overset{f}{\to} B \overset{g}{\to} M \overset{h}{\to}
\label{magic1}
\end{equation}
\begin{equation}
M \overset{f'}{\to} B' \overset{g'}{\to} M^{\ast}
\overset{h'}{\to} \label{magic2}
\end{equation}
be the exchange triangles corresponding to $M$ and $M^{\ast}$
(see~\cite[\S6]{bmrrt06}), so that $B\to M$ is a minimal right
$\text{add}(\overline{T}')$-approximation of $M$ in $\mathcal{C}$ and
$B' \to M^{\ast}$ is a minimal right
$\text{add}(\overline{T}')$-approximation of $M^{\ast}$ in $\mathcal{C}$.
The following definition appears to be important for
understanding the link between representation theory and denominators
of cluster variables:

\begin{defin} \label{d:excomp}
Let $N$ be an exceptional indecomposable object of $\C$. We say that
$N$ is \emph{compatible with an exchange pair $(M,M^{\ast})$} if
the following holds whenever $M \not \simeq \tau N \not \simeq M^{\ast}$:
\begin{multline*}
\dim\Hom_{\C}(N,M)+\dim\Hom_{\C}(N,M^{\ast}) \\
=  \max(\dim\Hom_{\C}(N,B),\dim\Hom_{\C}(N,B')) 
\end{multline*}
If $N$ is compatible with every exchange pair $(M,M^{\ast})$ in
$\C$ we call $N$ \emph{exchange compatible}.
\end{defin}

\noindent We will investigate when this condition is satisfied.
Note that the case when $M\simeq \tau N$ or $M^{\ast}\simeq \tau N$ is
covered by the following.

\begin{lem}\label{l:oldc3}
Let $N$ be an exceptional indecomposable object of $\mathcal{C}$ and
suppose that $M\simeq \tau N$ or $M^{\ast}\simeq \tau N$.
Then we have that \begin{multline*}
\dim\Hom_{\C}(N,M)+\dim\Hom_{\C}(N,M^{\ast}) = \\
\max(\dim\Hom_{\C}(N,B),\dim\Hom_{\C}(N,B')) +1 \end{multline*}
\end{lem}

\begin{proof}
Assume that $M^{\ast} \simeq \tau N$ (the other case is similar).
Then there are exchange triangles
$M \to B \to \tau N \to$ and $\tau N \to B' \to M \to$.
We have that $\Hom_{\C}(N,M) \simeq D\Ext^1_{\C}(M, \tau N)$.
This space is one-dimensional, since 
$(M, \tau N)$ is an exchange pair (by \cite[7.5]{bmrrt06}).
Since $N$ is exceptional, $\Hom_{\C}(N,\tau N) = 0$, and we
have $\Hom_{\C}(N, B \amalg B') \simeq D \Ext^1(B \amalg B', \tau N)= 0$, since $B \amalg B' \amalg \tau N$
can be completed to a cluster-tilting object. 
\end{proof}

\begin{prop}\label{p:c1toc2}
Suppose that neither $M$ nor $M^{\ast}$ is isomorphic to $\tau N$.
Then the following are equivalent:
\begin{itemize}
\item[(a)] We have that \begin{multline*} \dim\Hom_{\C}(N,M)+\dim\Hom_{\C}(N,M^{\ast}) \\ =
\max(\dim\Hom_{\C}(N,B),\dim\Hom_{\C}(N,B')). \end{multline*}
\item[(b)]
Either the sequence
\begin{equation}
0\to \Hom_{\C}(N,M^{\ast})\overset{a}{\to} \Hom_{\C}(N,B)\overset{b}{\to}
\Hom_{\C}(N,M)\to 0
\label{firstsequence}
\end{equation}
is exact, or the sequence
\begin{equation}
0\to \Hom_{\C}(N,M)\overset{c}{\to} \Hom_{\C}(N,B')\overset{d}{\to}
\Hom_{\C}(N,M^{\ast})\to 0
\label{secondsequence}
\end{equation}
is exact.
\end{itemize}
\end{prop}

\begin{proof}

Consider the (not necessarily exact) sequence:
$$\Hom_{\C}(N,M^{\ast})
\overset{a}{\to} \Hom_{\C}(N,B)
\overset{b}{\to} \Hom_{\C}(N,M)
\overset{c}{\to} \Hom_{\C}(N,B')
\overset{d}{\to} \Hom_{\C}(N,M^{\ast}).
$$
By exactness at the second term we have that $\Ker b = \Im a$.
Let $r=\dim \Hom_{\C}(N,M^{\ast})$, $s=\dim \Hom_{\C}(N,B)$,
$t=\dim\Hom_{\C}(N,M)$ and $u=\dim\Hom_{\C}(N,B')$.
Then we have that
$s=\dim\Ker b+\dim\Im b=\dim\Im a+\dim \Im b \leq r+t$.
Similarly, by exactness at the fourth term, we have that $u\leq t+r$.

Then the sequence~\eqref{firstsequence} is exact if and only
if $b$ is surjective and $a$ is injective, if and only
if $\dim\Im b=t$ and $\dim\Im a=r$.
By the above argument that $s\leq r+t$, we see that this is true if and
only if $s=r+t$. Similarly, the sequence~\eqref{secondsequence}
is exact if and only if $u=r+t$. Hence (b) holds for $N$ if and only if
$s=r+t$ or $u=r+t$. Since in any case $s\leq r+t$ and $u\leq r+t$, we
see that (b) holds if and only $r+t=\max(s,u)$ as required.
\end{proof}

Let $X,Y$ be indecomposable objects of $\C$. We can always represent
such objects by indecomposable $kQ$-modules or objects of the form
$P_i[1]$ for some $i$. Call these also $X$ and $Y$, by abuse of
notation. By~\cite[1.5]{bmrrt06},
$\Hom_{\C}(X,Y)=\Hom_{\D}(X,Y)\amalg \Hom_{\D}(X,FY)$.
We call an element of $\Hom_H(X,Y)$ an \emph{$H$-map} from $X$ to $Y$, and
an element of $\Hom_{\D}(X,FY)$ an \emph{$F$-map} from $X$ to $Y$. We note
that the composition of two $H$-maps is an $H$-map and the composition of
two $F$-maps is zero. The composition of an $H$-map and an $F$-map is an
$F$-map.

\begin{prop}
If $N$ is an indecomposable exceptional object with $\End_{\C}(N) \not \simeq k$, then there is 
an indecomposable exceptional object $N^{\ast}$ such that $(N,N^{\ast})$ is an exchange pair and $N$ is not compatible with
this exchange pair. 
\end{prop}

\begin{proof}
Let $E$ be the Bongartz complement of $N$, as defined in \cite{bong}. Then $N \amalg E$ is a tilting module
such that $\Hom_H(N,E)=0$ (see \cite{happel}).

Now consider $N \amalg E$ as a cluster-tilting object in $\C$, and consider
the exchange triangles $N \to B' \to N^{\ast}$ and 
$N^{\ast} \to B  \to N$, where $B$ and $B'$ are in $\add E$. 

Then the map $N \to B'$ is a sum of $F$-maps. Since by assumption 
$\End_{\C}(N) \not\simeq k$, we have a non-zero $F$-map
$N \overset{\epsilon}{\to} N$. The composition $N \overset{\epsilon}{\to} N 
\to B'$ is zero, as the composition of $F$-maps is zero.
So $\Hom_{\C}(N,\ )$ applied to $N \to B' \to N^{\ast}$
does not give a short exact sequence.

It is clear that $\Hom_{\C}(N,\ )$ applied to $N^{\ast} \to B \to N$ does not give a short exact
sequence, since the identity map $N \to N$ does not factor through $B \to N$.

We claim that $N$ is not compatible with $(N, N^{\ast})$. We have that
$\Ext^1_{\C}(N,\tau N)\simeq D\Hom_{\C}(N,N)\not\simeq k$, so
$(N, \tau N)$ is not an exchange pair by~\cite[7.5]{bmrrt06}. It follows that
$N^{\ast} \not \simeq \tau N$. We also have that $N \not \simeq \tau N$,
since $N$ is exceptional.
So we are done by combining the above with Proposition \ref{p:c1toc2}.
\end{proof}

We have the following immediate consequence.

\begin{cor} \label{c:endpropc}
If $N$ is an exchange compatible object of $\C$, then $\End_{\C}(N) \simeq k$.
\end{cor}

There is the following description of the indecomposable exceptional modules $N$
with $\End_{\C}(N) \not \simeq k$.
We know that such $N$ must be regular. In the tame case the indecomposable
exceptional
modules in a tube of rank $t$ are those of quasilength at most $t-1$.
In the wild case, for a component of the AR-quiver there is some $r$ with $0 \leq r \leq n-2$, 
where $n$ is the number of simple modules, such that
$N$ is exceptional if and only if $N$ has quasilength at most $r$.

Let $X$ be an indecomposable regular module. Then we denote by $\Delta_X$ the
wing of $X$, as defined in \cite[3.3]{rin}.

The proof of the next result, which is based on
the proof of \cite[1.6]{kern}, \cite[9.5]{kern2}, was pointed out to
us by Kerner. 

\begin{prop} \label{p:endnotk}
Assume that $H$ is hereditary of infinite type, and that $N$ is an
indecomposable exceptional regular $H$-module. Then
$\End_{\C}(N) \not \simeq k$ if and only if $N$ has maximal quasilength
amongst the indecomposable exceptional modules in its component.
\end{prop}

\begin{figure}
$$ \xy \xymatrix@C=8pt@R=7pt{
&&& \tau N' \ar[dr] && N' \ar[dr] \\
&& \tau^2 N \ar[ur] \ar[dr] && \tau N \ar[ur] \ar[dr] && N \\
&&& \ar[ur] && \ar[ur] \\
\\
\\
A \ar@{--}[uuurrr] &&&&&&&& B \ar@{--}[uuulll] }
\endxy
$$
\caption{The proof of Proposition~\ref{p:endnotk}}
\label{f:endnotk}
\end{figure}
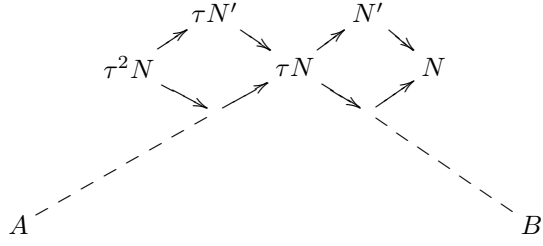

\begin{proof}
Let $N$ be indecomposable exceptional regular of quasilength $t$,
and let $N'$ be indecomposable with an irreducible epimorphism
$N\to N'$. Then $\End_{\C}(N) \not \simeq k$ if and only if
$\Hom_H(N, \tau^2 N) \neq 0$, and $N'$ is exceptional if and only
if $\Hom_H(N', \tau N') \simeq \Ext^1_H(N', N')= 0$, if and only
if $N$ does not have maximal quasilength amongst the indecomposable
exceptional modules in its component (note
that all indecomposable modules in a $\tau$-orbit of an AR-component of
a hereditary algebra are always either all exceptional or all not
exceptional).

Consider the exact sequences $0 \to A \to N' \to N \to 0$ and $0
\to \tau^2 N \to \tau N' \to B \to 0$. Since $N' \to N$ is an
epimorphism and $\tau^2 N \to \tau N'$ is a monomorphism, it
follows that $\Hom_H(N, \tau^2 N) \neq 0$ implies $\Hom_H(N' ,
\tau N') \neq 0$.

Assume conversely that there is a nonzero map $f \colon N' \to
\tau N'$. Since the wing $\Delta_{\tau N'}$ determined by $\tau
N'$ is standard \cite{kern2}, we have that $\Hom(A, \tau N') = 0$ (see
Figure~\ref{f:endnotk}). Hence the composition of the map $A\to
N'$ with $N\to \tau N'$ is zero. It follows that the map $N'\to
\tau N'$ factors through $N$ and therefore that there is a
non-zero map $f' \colon N \to \tau N'$.
Similarly, considering the wing $\Delta_{N'}$, we see that
$\Hom(N,B)=0$. Hence the composition of $f'$ with the map $\tau
N'\to B$ is zero, and $f'$ factors through $\tau^2 N$ and
therefore there is a non-zero map $N\to \tau^2 N$.
Hence we see that $\End_{\C}(N) \not \simeq k$ if and only if $N$
is exceptional of maximal quasilength.
\end{proof}

\begin{rem}
The proof in the tame case can be made even shorter.
\end{rem}

We note the following consequence of Lemma~\ref{c:endpropc} and
Proposition~\ref{p:endnotk}:

\begin{cor}\label{thm3.12}
Assume that $H$ is of infinite type.
If $N$ is indecomposable regular and has maximal quasilength amongst the
indecomposable exceptional modules in its component, then $N$ is not
exchange compatible.
\end{cor}

\section{Expressions for cluster variables}

In this section we show that for a given cluster-tilting object $T$ in the cluster category of $kQ$, all indecomposable direct summands 
$T_j$ of $T$ are exchange compatible if and only if all cluster variables in $\A(Q)$ have $T$-denominators.

Let $(\mathbf{x'},Q')$ be a seed, with $\mathbf{x}=\{x'_1,\ldots ,x'_n\}$, and assume
that each $x'_i$ has a $T$-denominator.
Let $T'_i=\alpha(x'_i)$ for
$i=1,2,\ldots ,n$ and set $T'=\amalg_{j=1}^n T'_j$.
By Theorem~\ref{bmrtsecond}, we have that $T'$ is a cluster-tilting object in $\C$
and that $Q'$ is the quiver of $\End_{\C}(T')^{\op}$.

Mutating $(\mathbf{x}',Q')$ at $x'_k$ we obtain a new cluster variable
$(x'_k)^{\ast}$.
Let $\overline{T}'=\amalg_{j\not=k}T'_j$.
Let $(T'_k)^{\ast}$ be the unique indecomposable object in $\C$ with
$(T'_k)^{\ast}\not\simeq T'_k$ such that $\overline{T}'\amalg (T'_k)^{\ast}$ is a
cluster-tilting object.

\begin{prop} \label{p:inductionstep}
With the above notation and assumptions,
the cluster variable $(x'_k)^{\ast}$ has a $T$-denominator if and only if
each summand $T_i$
of $T$ is compatible with the exchange pair $(T'_k,(T'_k)^{\ast})$.
\end{prop}

\begin{proof}
We consider the exchange triangles corresponding to the almost complete
cluster-tilting object $\overline{T}'$ in $\C$ (see triangles~\eqref{magic1}
and~\eqref{magic2}).
For an indecomposable direct summand $A$ of $T'$, we
denote the corresponding cluster variable in $\mathbf{x}'$ by $x_A$.
Note that all summands of $B$ and $B'$ are also summands
of $\overline{T}'$.
We therefore have, by assumption, that for any summand
$A$ of $B$ or $B'$ (or $A=T'_k$), either $A\not\simeq \tau T_i$ for any $i$,
and $x_A=f_A/t_A$, or $A\simeq \tau T_i$ for some $i$ and $x_A=f_Ay_i$,
where in either case $f_A$ is a polynomial satisfying the positivity
condition. We also define $M=T'_k$ and $M^{\ast}=(T'_k)^{\ast}$.
We note that by~\cite[3.1]{bmrt}, $\alpha((x'_k)^{\ast})=M^{\ast}$.

If $A_1,A_2,\ldots ,A_t$ are summands of $T'$, we write
$$x_{A_1\amalg A_2\amalg \cdots \amalg A_t}=x_{A_1}x_{A_2}\cdots x_{A_t}$$
and
$$f_{A_1\amalg A_2\amalg \cdots \amalg A_t}=f_{A_1}f_{A_2}\cdots f_{A_t}.$$
We also write
$$y_{\tau T_{i_1}\amalg \tau T_{i_2}\amalg \cdots \amalg
\tau T_{i_t}}=y_{i_1}y_{i_2}\cdots y_{i_t}$$
for any $i_1,i_2,\ldots ,i_t\in \{1,2,\ldots ,n\}.$

We have~\cite{bmrA} that
\begin{equation}
x_M(x_M)^{\ast}=x_B+x_{B'}, \label{mmstar}
\end{equation}
The argument now falls into different cases.

\noindent {\bf Case (I):} \\
Suppose first that neither $M$ nor $M^{\ast}$ is isomorphic to $\tau T_i$ for
any $i$. 
Write $B=B_0\amalg B_1$, where no summand of $B_0$ is of the form
$\tau T_i$ for any $i$.
Similarly, write $B'=B'_0\amalg B'_1$, where no summand of $B'_0$
is of the form $\tau T_i$ for any $i$.
Then by assumption $x_B=\frac{f_By_{B_1}}{t_{B_0}}$.
Similarly $x_{B'}=\frac{f_{B'}y_{B'_1}}{t_{B'_0}}$.

Let $m=\frac{\lcm(t_{B_0},t_{B'_0})}{t_{B_0}}$ and
$m'=\frac{\lcm(t_{B_0},t_{B'_0})}{t_{B'_0}}$. We have:
\begin{eqnarray*}
(x_M)^{\ast} & = & \frac{\frac{f_By_{B_1}}{t_{B_0}}+\frac{f_{B'}y_{B'_1}}{t_{B'_0}}}
{\frac{f_M}{t_M}} \\
& = & \frac{f_Bmy_{B_1}+f_{B'}m'y_{B'_1}}{\lcm(t_{B_0},t_{B'_0})f_M/t_M} \\
& = & \frac{(f_Bmy_{B_1}+f_{B'}m'y_{B'_1})/f_M}{\lcm(t_B,t_{B'})/t_M},
\end{eqnarray*}
noting that $t_B=t_{B_0}$ since $\Hom_{\C}(T_i,\tau T_j)=0$ for all $i,j$,
and similarly $t_{B'}=t_{B'_0}$.

Assume now that each of the $T_i$ is compatible
with $(M,M^{\ast})$, so that we have 
\begin{eqnarray*}
t_Mt_{M^{\ast}} & = & \prod y_i^{\dim \Hom_{\C}(T_i,M) + \dim \Hom_{\C}(T_i,M^{\ast})}   \\
 & = & \prod y_i^{\max(\dim\Hom_{\C}(N,B),\dim\Hom_{\C}(N,B'))}  = \lcm(t_B,t_{B'}). \end{eqnarray*}
Hence
\begin{equation}
(x_M)^{\ast}=\frac{(f_Bmy_{B_1}+f_{B'}m'y_{B'_1})/f_M}{t_{M^{\ast}}}.
\label{mstar}
\end{equation}
By definition of least common multiple, $m$ and $m'$ are coprime.
Since $B$ and $B'$ have no common factors~\cite[6.1]{bmrA},
$y_{B_1}$ and $y_{B'_1}$ are coprime.
Suppose that $m$ and $y_{B'_1}$ had a common factor $y_i$. Then there would
be a summand $X$ of $B'_0$ such that $\Hom_{\C}(T_i,X)\not=0$, and $\tau T_i$
was a summand of $B'$. But then
$$\Ext^1_{\C}(X,\tau T_i)\simeq D\Hom_{\C}(\tau T_i,\tau X)
\simeq D\Hom_{\C}(T_i,X)\not=0,$$ contradicting the fact that $B$ is the
direct sum of summands of a cluster-tilting object. Therefore $m$ and
$y_{B'_1}$ are coprime, and similarly $m'$ and $y_{B_1}$ are coprime. It follows
that $my_{B_1}$ and $m'y_{B'_1}$ are coprime.

Let $i\in \{1,2,\ldots ,n\}$. It follows from our assumptions that
$f_B(e_i)>0$ and $f_{B'}(e_i)>0$.
It is clear that $(my_{B_1})(e_i)\geq 0$ and $(m'y_{B'_1})(e_i)\geq 0$.
Since $my_{B_1}$ and $m'y_{B'_1}$ are coprime, these two numbers cannot be
simultaneously zero, so $(f_Bmy_{B_1}+f_{B'}m'y_{B'_1})(e_i)>0$. It follows
that $f_Bmy_{B_1}+f_{B'}m'y_{B'_1}$ satisfies the positivity condition.
By assumption, $f_M$ also satisfies the positivity condition.

Since $(x_M)^{\ast}$ is a Laurent polynomial in
$y_1,y_2,\ldots ,y_n$ by the Laurent phenomenon~\cite[3.1]{fominzelevinsky02},
$(f_Bmy_{B_1}+f_{B'}m'y_{B'_1})/f_M$ is also a Laurent polynomial.
Since it is defined at $e_i$ for all $i$, it must be a polynomial.
By the above, it satisfies the positivity condition.
It follows that $(x_M)^{\ast}$ has a $T$-denominator.

Now assume that at least one summand of $T$ is not compatible with
$(M,M^{\ast})$. Then
$t_Mt_{M^{\ast}}=c\lcm(t_B,t_{B'})$ for some nontrivial Laurent monomial $c$.
We then obtain that
\begin{eqnarray*}
x_{M^{\ast}} & = &
\frac{(f_Bmy_{B_1}+f_{B'}m'y_{B'_1})/f_M}{\lcm(t_B,t_{B'})/t_M} \\
&=&
c\cdot \frac{(f_Bmy_{B_1}+f_{B'}m'y_{B'_1})/f_M}{t_{M^{\ast}}},
\end{eqnarray*}
where, as before, $(f_Bmy_{B_1}+f_{B'}m'y_{B'_1})/f_M$ satisfies
the positivity condition.
It follows that the numerator
$c\cdot (f_Bmy_{B_1}+f_{B'}m'y_{B'_1})/f_M$ is either not a polynomial
or does not satisfy the positivity condition and thus that $x_{M^{\ast}}$ does
not have a $T$-denominator.

\noindent {\bf Case (II):} \\
Next, suppose that $M\simeq \tau T_i$ for some $i$.
Since $\Ext^1_{\C}(T_r,T_s)=0$ for all $r,s$, we cannot have $M^{\ast}\simeq \tau T_j$
for any $j$ (as $\Ext^1_{\C}(M,M^{\ast})\not=0$).
Assume that the $T_j$ are all compatible with the exchange pair $(M,M^{\ast})$. Then
we have
$t_{M^{\ast}}=y_i\lcm(t_B,t_{B'})$, noting that $t_M=1$. Here we use, in particular,
that 
\begin{multline*}
\dim\Hom_{\C}(N,M)+\dim\Hom_{\C}(N,M^{\ast}) \\ =
\max(\dim\Hom_{\C}(N,B),\dim\Hom_{\C}(N,B'))+1
\end{multline*}
holds for $N= T_i$, by Lemma \ref{l:oldc3}.
As in the above case, we have $x_M(x_M)^{\ast}=x_B+x_{B'}$. The same analysis,
together with the fact that $x_M=f_My_i$ for some polynomial $f$ in the
$y_j$ satisfying the positivity condition, provides:
\begin{eqnarray*}
(x_M)^{\ast} & = & \frac{(f_Bmy_{B_1}+f_{B'}m'y_{B'_1})/f_M}{\lcm(t_B,t_{B'})y_i} \\
& = & \frac{(f_Bmy_{B_1}+f_{B'}m'y_{B'_1})/f_M}{t_{M^{\ast}}}.
\end{eqnarray*}
As above, we see that the numerator is a polynomial which satisfies the
positivity condition, and it follows that $(x_M)^{\ast}$ has a $T$-denominator.

Now assume that at least one summand of $T$ is not compatible with
$(M,M^{\ast})$. Summands $T_k \not \simeq T_i$ which are not compatible with 
$(M,M^{\ast})$, will not satisfy $\dim\Hom_{\C}(N,M)+\dim\Hom_{\C}(N,M^{\ast})=
\max(\dim\Hom_{\C}(N,B),\dim\Hom_{\C}(N,B'))$ for $N= T_k$. 

It follows that $t_{M^{\ast}}=cy_i\lcm(t_B,t_{B'})$ for some nontrivial Laurent
monomial $c$.
Arguing as above we obtain that
\begin{eqnarray*}
x_{M^{\ast}}  &=& \frac{f_B m y_{B_1}+f_{B'} m' y_{B'_1}}{t_{M^{\ast}}/c} \\
&=& \left( \frac{f_B m y_{B_1}+f_{B'} m' y_{B'_1}}{t_{M^{\ast}}}\right) c.
\end{eqnarray*}
Since $f_B m y_{B_1}+f_{B'} m' y_{B'_1}$ satisfies the positivity
condition, it follows that the numerator
$(f_B m y_{B_1}+f_{B'} m' y_{B'_1}) c$
is either not a polynomial or does not satisfy the positivity condition
and therefore that $(x_{M})^{\ast}$ does not have a $T$-denominator.

\noindent {\bf Case (III):} \\
Next, suppose that $M^{\ast}\simeq \tau T_i$ for some $i$.
Since $\Ext^1_{\C}(T_r,T_s)=0$ for all $r,s$, we cannot have $M\simeq \tau T_j$
for any $j$ (as $\Ext^1_{\C}(M,M^{\ast})\not=0$).
Assume that the $T_j$ are all exchange compatible with $(M, M^{\ast})$. We have, using
Proposition \ref{p:c1toc2}, that
$t_{M}=y_i\lcm(t_B,t_{B'})$, noting that $t_{M^{\ast}}=1$.
As in the above case, we have $x_M(x_M)^{\ast}=x_B+x_{B'}$. The same analysis
provides:
\begin{eqnarray*}
(x_M)^{\ast} & = & \frac{(f_Bmy_{B_1}+f_{B'}m'y_{B'_1})/f_M}{\lcm(t_B,t_{B'})/t_M} \\
& = & \left( (f_Bmy_{B_1}+f_{B'}m'y_{B'_1})/f_M \right)y_i.
\end{eqnarray*}
As above, we see that $(f_Bmy_{B_1}+f_{B'}m'y_{B'_1})/f_M$
is a polynomial which satisfies the positivity condition and it follows that
$(x_M)^{\ast}$ has a $T$-denominator.

Now assume that at least one summand of $T$ is not compatible with
$(M,M^{\ast})$. Then, as in Case (II), it follows that
$t_M=cy_i\lcm(t_B,t_{B'})$ for some nontrivial Laurent
monomial $c$.
Arguing as above we obtain that
\begin{eqnarray*}
(x_M)^{\ast} & = & \frac{(f_Bmy_{B_1}+f_{B'}m'y_{B'_1})/f_M}{\lcm(t_B,t_{B'})/t_M} \\
& = & \left( (f_Bmy_{B_1}+f_{B'}m'y_{B'_1})/f_M \right)cy_i.
\end{eqnarray*}

Since $f_B m y_{B_1}+f_{B'} m' y_{B'_1}$ satisfies the positivity
condition, it follows that $(f_B m y_{B_1}+f_{B'} m' y_{B'_1}) c$
is either not a polynomial or does not satisfy the positivity condition
and therefore that $(x_{M})^{\ast}$ does not have a $T$-denominator.
\end{proof}

\begin{thm} \label{t:all}
All indecomposable direct summands $T_j$ of $T$ are exchange compatible
if and only if all cluster variables in $\mathcal{A}(Q)$ have a $T$-denominator.
\end{thm}

\begin{proof}
We note that $\alpha(y_i)=\tau T_i$ for $i=1,2,\ldots ,n$ and that the
quiver $\Gamma$ of the seed $(\mathbf{y},\Gamma)$ is the quiver of
$\End_{\C}(\amalg_{j=1}^n \tau T_j)^{\op}$ (see the start of
Section~\ref{s:mainresults}).
Suppose first that all indecomposable direct summands $T_j$ of $T$ are
exchange compatible. Let $x$ be a cluster variable of $\mathcal{A}(Q)$.
It follows from Proposition~\ref{p:inductionstep}
by induction on the smallest number of exchanges needed to
get from $\mathbf{y}$ to a cluster containing $x$ that $x$ has a $T$-denominator.

Suppose that there is an indecomposable direct summand of $T$ which
is not exchange compatible. 
Consider the cluster $\mathbf{z}$ of $\mathcal{A}(Q)$ which is
a minimal distance in the exchange graph from $\mathbf{y}$ such that 
some summand of $T$ is not compatible with an
exchange pair $(\alpha(x),\alpha(x^{\ast}))$ of modules coming from an
exchange pair $(x,x^{\ast})$ of cluster variables with $x\in\mathbf{z}$.
Arguing as above, we see that the cluster variables in $\mathbf{z}$
all have a $T$-denominator. By Proposition~\ref{p:inductionstep} we
see that $x^{\ast}$ does not have a $T$-denominator.
\end{proof}

The fact that $\alpha$ is a bijection gives us more information, since it
follows that $\alpha(x)=\tau T_i$ if and only if $x=y_i$. For the rest of
this article, we freely use this fact. We therefore have:

\begin{cor} \label{c:correctform}
All indecomposable direct summands $T_j$ of $T$ are exchange compatible
if and only if for every cluster variable $x$ of $\mathcal{A}(Q)$ we have: \\
(a) The object $\alpha(x)$ is not isomorphic to $\tau T_i$ for any $i$, and
$x=f/t_{\alpha(x)}$ for some polynomial $f$ in the $y_i$ satisfying the
positivity condition, or \\
(b) We have that $\alpha(x)\simeq \tau T_i$ for some $i$ and $x=y_i$.
\end{cor}

\begin{rem} \rm
Recall that the exponent of $y_i$ in $t_M$ is $d_i=\dim
\Hom_{\C}(T_i,M)$. It follows from~\cite[Theorem A]{bmrB} that
$d_i$ is the multiplicity of the simple top of the
projective $\L=\End_{\C}(T)^{\op}$-module $\Hom_{\C}(T,T_i)$
as a composition factor in $\Hom_{\C}(T,M)$: we have that
\begin{eqnarray*}
\Hom_{\C}(T_i,M) & \simeq & \Hom_{\C/\add(\tau T)}(T_i,M) \\
& \simeq & \Hom_{\L}(\Hom_{\C}(T,T_i),\Hom_{\C}(T,M)),
\end{eqnarray*}
using the fact that $\Hom_{\C}(T_i,\tau T)=0$ (as $T$ is a
cluster-tilting object in $\C$).
Thus $(d_1,d_2,\ldots ,d_n)$ is the dimension vector of the $\Lambda$-module
$\Hom_{\C}(T,M)$.
\end{rem}

\section{The preprojective case} \label{s:thepreprojectivecase}

In this section we adopt the notation from the previous section,
but assume in addition that $N$ is a preprojective module, a preinjective
module, or the shift of a projective module. Our aim is to show that
such modules are exchange compatible.

\begin{prop} \label{p:oneexact}
Suppose that neither $M$ nor $M^{\ast}$ is isomorphic to $\tau N$. Then
at least one of the following holds (where the maps are
induced from triangle~\eqref{magic1} or triangle~\eqref{magic2}
respectively): \\
(a) The sequence
$$0\to \Hom_{\C}(N,M^{\ast})\overset{a}{\to} \Hom_{\C}(N,B)\overset{b}{\to} \Hom_{\C}(N,M)\to 0$$
is exact, or \\
(b) The sequence
$$0\to \Hom_{\C}(N,M)\overset{c}{\to} \Hom_{\C}(N,B')\overset{d}{\to} \Hom_{\C}(N,M^{\ast})\to 0$$
is exact.
\end{prop}

\begin{proof}
Assume first that $N$ is projective.

\noindent {\bf Case (I):} \\
Suppose that $M$ and $M^{\ast}$ are
both modules. By~\cite{bmrrt06}, we know that either triangle~\eqref{magic1}
or triangle~\eqref{magic2}
arises from a short exact sequence of modules. Suppose first
that triangle~\eqref{magic1} arises from a short exact sequence
$$0\to M^{\ast}\to B \to M\to 0$$
of modules.
Applying the exact functor $\Hom_H(N,-)$ to this, we obtain the exact
sequence:
$$0\to \Hom_H(N,M^{\ast})\to \Hom_H(N,B)\to \Hom_H(N,M)\to 0,$$
noting that $\Ext^1_H(N,M^{\ast})=0$ as $N$ is projective.
Since $\Hom_{\C}(N,X)\simeq \Hom_H(N,X)$ for all modules $X$
(cf.~\cite[1.7(d)]{bmrrt06}), we obtain (a) above.
Similarly, if triangle~\eqref{magic2} arises from a short exact sequence
of modules, we obtain (b) above.

\noindent {\bf Case (II):} \\
Suppose that $M=P[1]$ is the shift of a projective indecomposable
module $P$.
Since $\Ext^1_{\C}(P[1],P'[1])=0$ for any projective module $P'$, we must
have that $M^{\ast}$ is a module.
Let $I$ be the indecomposable injective module corresponding to $P$.
Since $H$ is the path algebra of an acyclic quiver, either $\Hom_{\C}(N,P)=\Hom_H(N,P)=0$ or
$\Hom_{\C}(N,I)=\Hom_H(N,I)=0$. If $\Hom_{\C}(N,P)=0$, then,
applying $\Hom_{\C}(N,\ )$ to triangle~\eqref{magic1}, we obtain the
exact sequence
$$0\to \Hom_{\C}(N,M^{\ast})\to \Hom_{\C}(N,B)\to
0,$$
noting that $\Hom_{\C}(N,M)=\Hom_{\C}(N,P[1])=0$ as $N$ is projective.
We see that
(a) follows.

If $\Hom_{\C}(N,I)=0$, we note that $M[1]=P[2]\simeq \tau P[1]\simeq I$,
so applying $\Hom_{\C}(N,\ )$ to triangle~\eqref{magic2}, we obtain
the exact sequence
$$0\to \Hom_{\C}(N,B')\to \Hom_{\C}(N,M^{\ast})\to 0,$$
noting that $\Hom_{\C}(N,M)=\Hom_{\C}(N,P[1])=0$, and (b) follows.

The argument in case $M^{\ast}$ is not a module is similar.

Now assume that $N$ is a preprojective or preinjective module, or
$N$ is a shift of a projective.
Then $\tau^n N$ is projective for some $n\in\mathbb{Z}$.
Since $M\not\simeq \tau N$, then $\tau^n M\not\simeq \tau (\tau^n N)$, and
similarly for $M^{\ast}$.
By the above, either (a) or (b) holds for $\tau^n(N)$, so one of
$$0\to \Hom_{\C}(\tau^n N,\tau^n M^{\ast})\to \Hom_{\C}(\tau^n N,\tau^n B)
\to \Hom_{\C}(\tau^n N,\tau^n M)\to 0,$$
or
$$0\to \Hom_{\C}(\tau^n N,\tau^n M)\to \Hom_{\C}(\tau^n N,\tau^n B')
\to \Hom_{\C}(\tau^n N,\tau^n M^{\ast})\to 0$$
is exact. Hence one of
$$0\to \Hom_{\C}(N,M^{\ast})\to \Hom_{\C}(N,B)
\to \Hom_{\C}(N,M)\to 0,$$
or
$$0\to \Hom_{\C}(N,M)\to \Hom_{\C}(N,B')
\to \Hom_{\C}(N,M^{\ast})\to 0$$
is exact (since $\tau$ is an autoequivalence of $\C$) and we are done.
\end{proof}

\begin{cor} \label{c:preprojectiveC}
Each indecomposable preprojective module, preinjective module,
or shift of an indecomposable projective module is exchange compatible.
\end{cor}

\begin{proof} This follows from Proposition~\ref{p:oneexact}. 
\end{proof}

\vskip 0.2in
\noindent {\bf Proof of Theorem~\ref{main1}:} We observe that
Theorem~\ref{main1}(a) follows from Corollary~\ref{c:preprojectiveC} and
Theorem~\ref{main1}(b) follows from Corollary~\ref{c:endpropc}.~$\Box$

\vskip 0.2in
\noindent {\bf Proof of Corollary~\ref{main3}:} 
If $Q$ is Dynkin, then $kQ$ has no regular modules, and if $Q$ has only two
vertices, there are no exceptional regular modules. It follows from
Theorem~\ref{main1}(a) that in these cases, every cluster variable has
a $T$-denominator for every cluster-tilting object $T$.
Conversely, suppose that every cluster variable has a $T$-denominator for every cluster-tilting object $T$.
We note that
any exceptional indecomposable $kQ$-module can be completed to a tilting module
and therefore to a cluster-tilting object by~\cite[3.3]{bmrrt06}.
If $Q$ was not Dynkin and had more than two
vertices, $kQ$ would have an indecomposable regular exceptional module $N$ of
maximal quasilength amongst the indecomposable exceptional modules in its
AR-quiver component (see \cite{rin2} for the wild case).
By Proposition~\ref{p:endnotk}, we would have $\End_{\C}(N)\not\simeq k$, a
contradiction to Theorem~\ref{main1}(b).~$\Box$

\section{Criteria for exactness}

In this section we shall interpret exchange compatibility further.
The results here hold for any finite dimensional path algebra. They will
be used in the tame case in the next section.
To ease notation we let $\Hom_{\D}(\ , \ ) = (\ , \ )$ and
$\Ext_{\D}^1(\ , \ ) = {}^1(\ ,\ )$ (note that these coincide with the
corresponding $\Hom$ and $\Ext^1$-spaces over $H$ if both objects are
$H$-modules). We also let
$\Hom_{\C_H}(\ , \ ) = (\ , \ )_{\C}$ and
$\Ext^1_{\C_H}(\ , \ ) = {}^1(\ , \ )_{\C}$.

\begin{lem}\label{lem3.3}
Let $0 \to M^{\ast} \overset{f}{\to} B \overset{g}{\to} M \to 0$ be an exact
exchange sequence of
$H$-modules, and let $N$ be an exceptional indecomposable object in $\C$.
\begin{itemize}
\item[(a)]{$(N,\ )_{\C}$ gives a right exact sequence if and only if any $H$-map $h
\colon N \to M$
factors through $g$.}
\item[(b)]{$(N,\ )_{\C}$ gives a left exact sequence if and only if any $H$-map $s
\colon \tau N \to M$
factors through $g$.}
\end{itemize}
\end{lem}

\begin{proof}
Applying $(N,\ )_{\C}$ gives
the long exact sequence

$$
\xy
\xymatrix@C=9pt@R=3pt{
& & & 0 \ar[r] & (N,M^{\ast}) \ar[r] & (N,B) \ar[r] & (N,M) \ar[r] & \\
& & & & \amalg & \amalg & \amalg & \\
0 \ar[r] & (\tau N,M^{\ast}) \ar[r] & (\tau N,B) \ar[r] & (\tau N,M) \ar[r] &
{}^1(\tau N,M^{\ast}) \ar[r] & {}^1(\tau N,B) \ar[r] & {}^1(\tau N,M) \ar[r] & 0
}
\endxy
$$
The claim follows directly from this.
\end{proof}

\begin{rem} \rm \label{assumepreproj}
Given a finite number of indecomposable objects in a cluster category, where each of them is
preprojective, preinjective or the the shift of a projective, then
by changing the hereditary algebra $H$,
if necessary (i.e.\ via a tilt), we can assume that in fact all objects are 
preprojective (or preinjective). 
We make use of this in the sequel.
\end{rem}

\begin{lem}\label{lem3.4}
Let $M \to B' \to M^{\ast} \to$ be an exchange triangle, where $M$
is regular and $M^{\ast}$ is preprojective. Let $N$ be an indecomposable regular
exceptional module.
\begin{itemize}
\item[(a)]{$(N,\ )_{\C}$ gives a left exact sequence if and only if for any nonzero
$F$-map $N \to M$ the
composition $N \to M \to B'$ is nonzero.}
\item[(b)]{$(N,\ )_{\C}$ gives a right exact sequence if and only if for any nonzero
$F$-map
$\tau^{-1}N \to M$ the
composition $\tau^{-1} N \to M \to B'$ is nonzero.}
\end{itemize}
\end{lem}

\begin{proof}
\noindent (a)
One direction is obvious. Assume that for any nonzero
$F$-map $N \to M$ the
composition $N \to M \to B'$ is nonzero, 
then we need to see that the map $(N,M)_{\C}\to (N,B')_{\C}$ is injective.
Any non-zero map $N\to M$ in $\C$ must have the form $(h,h')$, where $h$
is an $H$-map and $h'$ is an $F$-map, by~\cite[1.5]{bmrrt06}.
Suppose, for a contradiction, that its composition with the map from
$M$ to $B'$ is zero. Then there is a map $s:N\rightarrow \tau^{-1}M^{\ast}$
such that $(h,h')=ts$, where $t:\tau^{-1}M^{\ast}\rightarrow M$.
Since $\tau^{-1}M^{\ast}$ is preprojective and
$M$ is regular, we have that $(\tau^{-1}M^{\ast},FM)=(\tau^{-1}M^{\ast},\tau^{-1}M[1])\simeq
(M^{\ast},M[1])\simeq {}^1(M^{\ast},M)\simeq D(M,\tau M^{\ast})=0$ since $\tau M^{\ast}$ is
preprojective. It follows that $t$ is an $H$-map. Similarly,
$(N,\tau^{-1}M^{\ast})=0$ since $N$ is regular and $\tau^{-1}M^{\ast}$ is
preprojective, so $s$ is an $F$-map. Comparing $H$-maps and $F$-maps,
we see that $h=0$ and $h'=ts$. Hence the $F$-map $h'$ is non-zero while
its composition with the map from $M$ to $B'$ is zero, contradicting
the assumption in (a).

$$
\xy
\xymatrix{
& N \ar^{(h,h')}[d] \ar@{--}_{s}[dl] &  & \\
\tau^{-1}M^{\ast} \ar^{t}[r] & M \ar[r] & B'\ \ar[r] & M^{\ast}
}
\endxy
$$

\noindent (b)
There is a long exact sequence
$$(N,M)_{\C} \to (N,B')_{\C} \to (N,M^{\ast})_{\C} \to (N,\tau M)_{\C} \to (N,\tau
B')_{\C} \to
(N,\tau M^{\ast})_{\C}.$$
Hence $(N,B')_{\C} \to (N,M^{\ast})_{\C}$ is an epimorphism if and only if
$(N,\tau M)_{\C} \to (N,\tau B')_{\C}$ is a monomorphism. This clearly holds
if and only if $(\tau^{-1}N,M)_{\C} \to (\tau^{-1}N,B')_{\C}$ is a monomorphism.
By (a), this holds if and only if for any non-zero $F$-map $\tau^{-1}N \to M$,
the composition $\tau^{-1}N \to M \to B'$ is nonzero.
\end{proof}

We have the following direct consequence of Lemmas~\ref{lem3.3}
and~\ref{lem3.4}.

\begin{prop}\label{prop3.7}
Let $(M, M^{\ast})$ be an exchange pair with $M$ regular and $M^{\ast}$
preprojective.
Let $0 \to M^{\ast} \to B \to M \to 0$ and $M \to B' \to M^{\ast} \to$
be the corresponding exchange sequence and triangle, where
$B = B_1 \amalg M_1$ and $B' = B'_1 \amalg M_1'$ with $B_1, B_1'$ preprojective and
$M_1, M_1'$ regular.
Let $N$ be an indecomposable regular exceptional module.
\begin{itemize}
\item[(a)]{Applying $(N,\ )_{\C}$ to $0 \to M^{\ast} \to B_1 \amalg M_1 \to M \to 0$
gives a short exact sequence
if and only if
any $H$-map $N \amalg \tau N \to M$ factors through $M_1 \to M$.}
\item[(b)]{Applying $(N,\ )_{\C}$ to $M \to B'_1 \amalg M'_1 \to M^{\ast} \to $
gives a short exact sequence if and only if
for any nonzero $F$-map $N \amalg \tau^{-1} N \to M$ the composition $N \amalg
\tau^{-1} N \to M \to M'_1$ is nonzero.
}
\end{itemize}
\end{prop}

\begin{proof}
\noindent (a)
If $N \amalg \tau N \to M$ is an $H$-map, it cannot factor through $B_1 \to M$. We apply Lemma~\ref{lem3.3}.

\noindent (b)
Since $B'_1$ is preprojective, the map $M \to B'_1$ is an $F$-map. Hence if $N \amalg \tau^{-1} N \to M$ 
is an $F$-map, then the
composition $N \amalg \tau^{-1} N \to M \to B'_1$, must be zero.
We apply Lemma~\ref{lem3.4}.
\end{proof}

\section{The tame case} \label{s:thetamecase}

In this section, we investigate the case in which $H$ is tame in more detail. 
In particular, we determine which exceptional regular modules are exchange 
compatible in this case.

\begin{lem}\label{lem3.6}
Let $H$ be a tame hereditary algebra and let $(M, M^{\ast})$ be an exchange pair with $M$ regular and $M^{\ast}$
preprojective.
Let $0 \to M^{\ast} \to B \to M \to 0$ and $M \to B' \to M^{\ast} \to$
be the corresponding exchange sequence and triangle, where
$B = B_1 \amalg M_1$ and $B' = B'_1 \amalg M_1'$ with $B_1, B_1'$ preprojective and
$M_1, M_1'$ regular.
Then the maps $M_1 \to M$ and $M \to M_1'$ are $H$-maps.
\end{lem}

\begin{proof}
The first claim holds obviously, since
the exchange sequence is in $\mod H$. For the second claim, by
Remark~\ref{assumepreproj}, by replacing $H$ with $H'$ via a tilt,
we can assume that $M^{\ast}$ is preinjective,
and similarly that $B'_1$ is preinjective.
Then the triangle $M \to B' \to M^{\ast} \to$
is the image of an exact sequence of $H'$-modules
$0 \to {}_{H'} M \to {}_{H'} B' \to {}_{H'} M^{\ast} \to 0$.
Then $M \to M'_1$ is an $H'$-map and consequently also an $H$-map.
\end{proof}

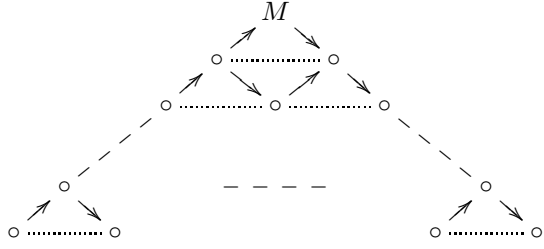
\begin{figure}
$$
\xy
\xymatrix@C=8pt@R=7pt{
&&&&& M \ar[dr] &&&& \\
&&&& \circ \ar[ur] \ar[dr] \ar@{.}[rr] && \circ \ar[dr] &&& \\
&&& \circ \ar@{--}[ddll] \ar[ur] \ar@{.}[rr] && \circ \ar[ur] \ar@{.}[rr] && \circ
\ar@{--}[ddrr] && \\
&&&&&&&&& \\
&\circ \ar[dr] &&& \ar@{--}[rr] &&&&& \circ \ar[dr] & \\
\circ \ar[ur] \ar@{.}[rr] && \circ &&&&&& \circ
     \ar[ur] \ar@{.}[rr] && \circ
}
\endxy
$$
\caption{The AR-quiver of the subcategory of subfactors of $M$}
\label{f:Matthetop}
\end{figure}

\begin{rem} \rm \label{Ansubcategory} \rm
Let $\E$ denote a category which is either the module category
of the path algebra $\L = \L_n$ of a quiver of type $A_n$ with linear orientation, or the abelian category
associated with a
tube, i.e.\ a connected component
of regular modules over some tame hereditary algebra $\L$.
For an indecomposable exceptional object $M$ in $\E$, we let $\Delta_M$ denote the
additive subcategory generated by the
subfactors of $M$ inside $\E$, or equivalently
the extension closure of $\Sub M \cup \Fac M$. The indecomposable objects
in $\Delta_M$ form a full subquiver of the AR-quiver of $\L$, which is a
triangle with $M$ sitting on the top (see Figure~\ref{f:Matthetop}).
\end{rem}

\begin{lem}\label{lem3.8}
Consider the quiver of type $A_t$ with linear orientation and let $\L$ be its
path algebra.
Let $B$ be the indecomposable projective and injective $\L$-module.
Assume
$B_1 \to B_2 \to \cdots \to B_s \to B \to C_1 \to \cdots \to C_{t-s}$ is a
sequence of irreducible maps, where $t>s>0$.
If $T$ is a tilting module over $\L$, it must contain one of the $B_i$
(for $1 \leq i \leq s$) or $C_j$ (for $1 \leq j \leq t-s$) as an
indecomposable direct summand.
\end{lem}

\begin{proof}
Assume to the contrary that $T$ is a tilting module having none
of the $B_i, C_j$ as indecomposable direct summands.
The object $B$ is necessarily a summand of $T$, since it is projective-injective.
Since $B_s$ only has extensions with the the proper factors of $B$, at
least one such factor must be a summand of $T$.
Let $E_j \in \Fac B \cap \add (T/B)$,
with length $j$ maximal. Then, by assumption $j \leq t-(t-s)-1 = s-1$.

Let $P_i$ be the indecomposable projective of length $i$,
and let $k=t-j-1$. Then $T/B$ is in $\Delta_{E_j} \cup \Delta_{P_k}$,
since, for $k <l <t$, $P_l$ and all factors of $P_l$ which are not in $\Delta_{E_j}$ 
have non-trivial extensions with $E_j$.
Note that $T/B$ has $t-1$ nonisomorphic indecomposable summands.
The subcategory $\Delta_{E_j}$ has clearly at most $j$ summands in $T$,
so at least $t-1-j=k$ summands $T_1, \dots, T_{t-(j+1)}$of $T$ are in
$\Delta_{P_k}$.
Note that $\Delta_{P_k}$ is clearly equivalent to $\mod \L_k$, so
$T'= T_1 \amalg \dots \amalg T_{t-(j+1)}$ forms a tilting module in this subcategory,
and hence $P_k$ is a summand in $T'$. But this is a contradiction, since
$P_k$ is equal to some
$B_i$.
We illustrate the situation in Figure~\ref{f:typeA}.

\begin{figure}
\begin{center}
\includegraphics{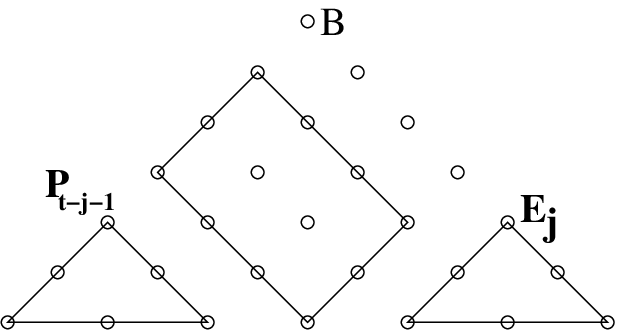} 
\end{center}
\caption{Proof of Lemma~\ref{lem3.8}.}
\label{f:typeA}
\end{figure}

\end{proof}

\begin{rem} \rm \label{Ancover}
Let $M$ be an indecomposable regular module in a tube $\T$ of rank $n$
for a tame hereditary algebra $H$.
Assume $\ql M = t \leq n-1 = \rank \T$ (where $\ql$ denotes quasilength).
Let $\L_t$ be the
path algebra of a quiver of type $A_t$ with linear orientation.
Then there is a functor
$F_t \colon \mod \L_t \to \Delta_M \subseteq \T$, which is clearly fully faithful
(from the structure of the tube), exact
and dense and preserves indecomposability.
\end{rem}

The following is due to Strauss \cite[Cor. 3.7]{strauss}.
\begin{lem}\label{lem3.9}
Let $M$ be an indecomposable exceptional $H$-module
for $H$ tame hereditary, lying in a tube $\T$ of rank $n$.
Let $\ql(M) = t$ (note that necessarily $t\leq n-1$).
Let $T = T' \amalg M$ be a tilting $H$-module.
Then $\add T \cap \Delta_M = \add \widetilde{T_0}$, where
$\widetilde{T_0} = F_t(T_0)$,
for a tilting module $T_0$ in $\mod \L_t$.
\end{lem}

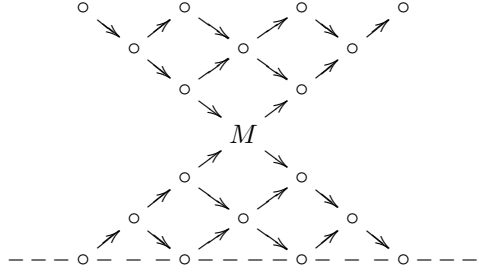
\begin{figure}
$$
\xy
\xymatrix@C=8pt@R=5pt{
&&\circ \ar[dr] && \circ \ar[dr] && \circ \ar[dr] && \circ \\
&&&\circ \ar[ur] \ar[dr] && \circ \ar[ur] \ar[dr] && \circ \ar[ur] \\
&&&&\circ \ar[ur] \ar[dr] && \circ \ar[ur] \\
&&&&&M \ar[ur] \ar[dr] \\
&&&& \circ \ar[dr] \ar[ur] && \circ \ar[dr] \\
&&& \circ \ar[dr] \ar[ur] && \circ \ar[dr] \ar[ur] && \circ \ar[dr] \\
&& \ar@{--}[ll] \ar@{--}[rr] \circ \ar[ur] && \circ \ar@{--}[rr] \ar[ur] && \circ
\ar@{--}[rr] \ar[ur] && \circ \ar@{--}[rr] &&
}
\endxy
$$
\caption{Proof of Lemma~\ref{lem3.9}. The dotted line indicates the
$\tau$-orbit of quasi-simple objects in the tube.}
\label{f:tubeorthogonal}
\end{figure}

\begin{lem}\label{lem3.10}
Let $\T$ be a tube of rank $n\geq 3$. Let $N,M$ be indecomposable objects
in $\T$ with $\ql(N)=j$, where $1\leq j\leq n-2$, and $\ql(M)=t$, where
$1\leq t\leq n-1$. Suppose that
$$X_1\to X_2\to \cdots \to X_t=M$$
is a sequence of irreducible maps with $X_t=M$ and $\ql(X_i)=i$
and that $$M=Y_t\to Y_{t-1}\to \cdots \to Y_1$$
is a sequence of irreducible maps with $Y_t=M$ and $\ql(Y_i)=i$.

(a) Let $f\in\Hom_{\T}(N,M)$ and $g\in\Hom_{\T}(M,\tau^2 N)$ be non-zero maps.
Then the number of modules $X_i$, $1\leq i<t$, that $f$ factors through plus
the number of modules $Y_i$, $1\leq i<t$, that $g$ factors through is at least
$t$. \\
(b) Statement (a) holds if we assume instead that $f\in\Hom_{\T}(\tau N,M)$
or $g\in \Hom_{\T}(M,\tau N)$ or both.
\end{lem}

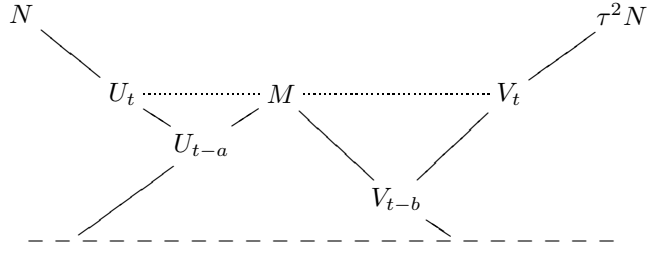
\begin{figure}
$$
\xy
\xymatrix@C=8pt@R=5pt{
N \ar@{-}[ddrr] & & & & & & & & & & \tau^2 N \ar@{-}[ddll] \\
& & & & & & & & & & \\
& & U_t \ar@{-}[dr] \ar@{.}[rr] & & M \ar@{-}[dl] \ar@{-}[ddrr] \ar@{.}[rrrr] & &
& & V_t \ar@{-}[ddll] \\
& & & U_{t-a} \ar@{-}[ddll] \\
& & & & & & V_{t-b} \ar@{-}[dr] \\
\ar@{--}[rrrrrrrrrr] & & & & & & & &  & &
}
\endxy
$$
\caption{Proof of Lemma~\ref{lem3.10} (the dotted line indicates the
$\tau$-orbit of quasi-simple objects in the tube).}
\label{f:tubefactoring}
\end{figure}

\begin{proof}
(a) Suppose first that $j\geq t$. Let
$$N=U_j\to U_{j-1}\to \cdots \to U_1$$
be a sequence of irreducible maps with $N=U_j$ and $\ql(U_i)=i$ for all $i$.
Then $U_t=\tau^a(M)$ for some $a$, where $1\leq a\leq n-1$,
so $U_{t-a}=X_{t-a}$. From the
structure of the tube it follows that $f$ factors through $X_{t-a},
X_{t-a+1},\ldots ,X_{t-1}$, a total of $a$ modules.

Let
$$V_1\to V_2\to \cdots \to V_j=\tau^2 N$$
be a sequence of irreducible maps with $V_j=\tau^2 N$ and
$\ql(V_i)=i$ for all $i$.
Then $V_t=\tau^{-b}(M)$ for some $b$, where $1\leq b\leq n-1$,
so $V_{t-b}=Y_{t-b}$. From the
structure of the tube it follows that $g$ factors through $Y_{t-b},
Y_{t-b+1},\ldots ,Y_{t-1}$, a total of $b$ modules.

Since $\tau^{n-2}(\tau^2(N))=N$ and $\ql(N)=j$, we have that
$\tau^{n-2-(j-t)}(V_t)=U_t$, so $a+b=n-2-(j-t)=(n-2-j)+t\geq t$
since we have assumed $j\leq n-2$. See Figure~\ref{f:tubefactoring}.

The argument in case $j<t$ is similar, so we are done.

For (b) we just need to remark that a non-zero map $\tau N \to M$
factors through more of the $X_i$ than $f$ does, and that a non-zero
map $M\to \tau N$ factors through more of the $Y_i$ than $g$
does. Therefore (a) holds if $N$ is replaced by $\tau N$, or $\tau^2 N$
is replaced by $\tau N$, or both. Statement (b) follows.
\end{proof}

\begin{lem}\label{lem3.11}
Let $H$ be tame and $N$ indecomposable exceptional in a tube of rank $n \geq 2$, so
that
$\ql (N) \leq n-1$. Then $\ql (N) = n-1$ if and only if
$\dim_k \End_{\C}(N)>1$. In
this case, there is a non-zero $F$-map from $N$ to $N$ and  $\dim_k \End_{\C}(N) =2$.
\end{lem}

\begin{proof}
It is easy to see that
$\Hom_H(N, \tau^2 N) \neq 0$ (and is one-dimensional)
if and only if $\ql(N) = n-1$, using
$\tau^2 N = \tau^{-(n-2)}N$.
Furthermore, $\Hom_{\C}(N,N) = \Hom(N,N) \amalg \Hom(N,FN)$.
We have that $\Hom(N,FN)\simeq \Ext^1(N, \tau^{-1}N)$,
and $\Ext^1(N, \tau^{-1}N) \simeq D\Hom(N, \tau^2 N)$.
\end{proof}

\begin{lem}\label{helpthm3.13}
Let $H$ be tame and $N$ indecomposable exceptional in a tube of rank
$n\geq 2$, with $\ql(N) = j \leq n-2$.
Let $(M,M^{\ast})$ be an exchange pair, where $M$ and $M^{\ast}$ are complements of some almost
complete cluster-tilting object $\overline{T}'$. We assume that $M$ is regular.
Then we have an exchange sequence $0 \to M^{\ast} \to B \to M \to 0$
and exchange triangle $M \to B' \to M^{\ast} \to $.
Write $B=B_1\amalg M_1$ and $B'=B'_1 \amalg M'_1$,
where $M_1,M_1'$ are regular and $B_1,B_1'$ are preprojective.
Then either:
\begin{itemize}
\item[(i)]{any $H$-map $N \amalg \tau N \to M$ factors through $B \to M$, or}
\item[(ii)]{for any nonzero $F$-map $N \amalg \tau^{-1}N \to M$, the composition
$N \amalg \tau^{-1}N \to M \to M'_1$ is nonzero.}
\end{itemize}
\end{lem}

\begin{proof}
We can clearly assume that there is some nonzero $H$-map
$N \amalg \tau N \to M$ and some nonzero $F$-map
$N \amalg \tau^{-1}N \to M$, and hence there is some nonzero
$H$-map $M \to \tau^2 N \amalg \tau N$. In particular,
$N$ and $M$ are in the same tube.
By Lemma~\ref{lem3.9}, $\overline{T}'$ is the direct sum of the
image under $F_t$ of a tilting module $T_0$ in $\mod \L_t$
and another summand (using the notation of Remark~\ref{Ancover}).
Let $\widetilde{T_0}=F_t(T_0)$.

\noindent {\bf Case (I): } \\
We first assume that there is a non-zero $H$-map from $N$ to $M$ and
a non-zero $H$-map from $M$ to $\tau^2 N$.
By Lemma~\ref{lem3.8}, combined with Lemma~\ref{lem3.10}(a), we see
that either some non-zero $H$-map from $N$ to $M$ factors through
some indecomposable direct
summand of $\overline{T}'$ or that some non-zero $H$-map
$M\to \tau^2 N$ factors through some indecomposable direct summand of
$\overline{T}'$.

In the first case, we note that, since $\Hom_H(N,M)\simeq k$, any non-zero
$H$-map from $N$ to $M$ factors through $\overline{T}'$.
Since any non-zero $H$-map from $\tau N$ to $M$ must factor through an
$H$-map from $N$ to $M$, it follows that any non-zero $H$-map from
$N\amalg \tau N$ to $M$ must factor through $\overline{T}'$, and therefore
through $B \to M$, and we see that (i) holds.

In the second case, arguing similarly, we see that any non-zero
$H$-map from $M$ to $\tau^2 N\amalg \tau N$ factors through the
$H$-map component of $M\to B'$ (noting that an $H$-map between two
modules cannot factor through an $F$-map).
It follows that the map from $\Hom(B',\tau^2 N\amalg \tau N)$ to
$\Hom(M,\tau^2 N\amalg \tau N)$ given by composition with the $H$-map
component of $M\to B'$ is surjective, so that the induced map from
$D\Hom(M,\tau^2 N\amalg \tau N))$ to
$D\Hom(B',\tau^2 N\amalg \tau N))$ is injective.
Hence the induced map from
$D\Hom(\tau^{-1} M,\tau N\amalg N))$ to
$D\Hom(\tau^{-1} B',\tau N\amalg N))$ is injective.
So the map from
$\Ext^1(N\amalg \tau^{-1}N,\tau^{-1}M)$ to
$\Ext^1(N\amalg \tau^{-1}N,\tau^{-1}B')$ is injective,
and therefore the map from
$\Hom(N\amalg \tau^{-1}N,\tau^{-1}M[1])$ to
$\Hom(N\amalg \tau^{-1}N,\tau^{-1}B'[1])$ is injective.
Hence the map from
$\Hom(N\amalg \tau^{-1}N,FM)$ to
$\Hom(N\amalg \tau^{-1}N,FB')$ is injective,
and we see that the composition of any non-zero $F$-map from
$N\amalg \tau^{-1}N$ to $M$ with the $H$-map component of
the map $M\to B'$ is non-zero. Note that composition of any
such $F$-map with the $F$-map component is always zero.
It follows that the
composition of any non-zero $F$-map from
$N\amalg \tau^{-1}N$ to $M$ with the map $M\to B'$ is
non-zero. We see that (ii) holds and we are done.

\noindent {\bf Case (II): } \\
We now assume that $\Hom_H(N,M)=0$, that there is a non-zero $H$-map from
$\tau N$ to $M$ and a non-zero $H$-map from $M$ to $\tau^2 N$.
By Lemma~\ref{lem3.8}, combined with Lemma~\ref{lem3.10}(b), we see
that either some non-zero $H$-map from $\tau N$ to $M$ factors through
some summand of $\overline{T}'$ or that some
non-zero $H$-map $M\to \tau^2 N$ factors through some summand of
$\overline{T}'$.

In the first case, we note that, since $\Hom(\tau N,M)\simeq k$, any non-zero
$H$-map from $\tau N$ to $M$ factors through $\overline{T}'$.
Since $\Hom_H(N,M)=0$, we see that any non-zero $H$-map from
$N\amalg \tau N$ to $M$ must factor through $\overline{T}'$, and therefore
through $B \to M$, and we see that (i) holds.

The argument in the second case is as in Case (I).

\noindent {\bf Case (III): } \\
We now assume that there is a non-zero $H$-map from $N$ to $M$,
that $\Hom_H(M,\tau^2 N)=0$ and that there is a non-zero $H$-map from
$M$ to $\tau N$.
By Lemma~\ref{lem3.8}, combined with Lemma~\ref{lem3.10}(b), we see
that either some non-zero $H$-map from $N$ to $M$ factors through
some summand of $\overline{T}'$ or that some
non-zero $H$-map $M\to \tau N$ factors through some summand of
$\overline{T}'$.

The argument in the first case is as in Case (I).

In the second case, since $\Hom_H(M,\tau^2 N)=0$, we see that any non-zero
$H$-map from $M$ to $\tau^2 N \amalg \tau N$ factors through the
$H$-map component of $M\to B'$ (noting that an $H$-map between two
modules cannot factor through an $F$-map). As in Case (I), we see
that (ii) holds.

\noindent {\bf Case (IV): } \\
Finally, we assume that $\Hom_H(N,M)=0$, $\Hom_H(M,\tau^2 N)=0$, that
there is a non-zero $H$-map from $\tau N$ to $M$ and that
there is a non-zero $H$-map from $M$ to $\tau N$.
By Lemma~\ref{lem3.8}, combined with Lemma~\ref{lem3.10}(b), we see
that either some non-zero $H$-map from $\tau N$ to $M$ factors through
some summand of $\overline{T}'$ or that some non-zero
$H$-map $M\to \tau N$ factors through some summand
of $\overline{T}'$.

In the first case we argue as in Case (II), and in the second case
we argue as in Case (III), and we are done.
\end{proof}

\begin{thm}\label{thm3.13}
Let $H$ be tame and $N$ indecomposable exceptional in a tube of
rank $n\geq 2$, with $\ql(N) = j \leq n-2$. Then $N$ is exchange compatible.
\end{thm}

\begin{proof}
Let $(M,M^{\ast})$ be an exchange pair, with $M$ and $M^{\ast}$ complements of some almost
complete cluster-tilting object $\overline{T}'$. 

\noindent {\bf Case (I): } \\
We assume first that $M$ is regular and $M^{\ast}$ preprojective.
Then we have an exchange sequence $0 \to M^{\ast} \to B \to M \to 0$
and exchange triangle $M \to B' \to M^{\ast} \to $.
Write $B=B_1\amalg M_1$ and $B'=B'_1 \amalg M'_1$,
where $M_1,M_1'$ are regular and $B_1,B_1'$ are preprojective.
By Lemma~\ref{helpthm3.13}, either
\begin{itemize}
\item[(i)]{any $H$-map $N \amalg \tau N \to M$ factors through $B \to M$, or}
\item[(ii)]{for any nonzero $F$-map $N \amalg \tau^{-1}N \to M$, the composition
$N \amalg \tau^{-1}N \to M \to M'_1$ is nonzero.}
\end{itemize}

If (i) holds, we see by Lemma~\ref{lem3.3} that the
sequence~\eqref{firstsequence} is exact.
If (ii) holds, we see by Lemma~\ref{lem3.4}
that the sequence~\eqref{secondsequence} is exact.

\noindent {\bf Case (II):} \\
Next, we assume that $M$ and $M^{\ast}$ are both preprojective.
We can still assume that
we have an exchange sequence $0 \to M^{\ast} \to B  \to M \to 0$
of $H$-modules, so that $B$ is also preprojective.
Applying $(N,\ )_{\C}$, arguing as in Lemma~\ref{lem3.3} and
using that $\Hom_H(N,X)= 0$ for any preprojective $X$
we see that the long exact sequence
$$(N, \tau^{-1}M)_{\C} \to (N,M^{\ast})_{\C} \to (N,B)_{\C} \to (N,M)_{\C} \to
(N,\tau M^{\ast})_{\C} \to$$
reduces to
$$0 \to \Ext^1_H(N, \tau^{-1}M^{\ast}) \to \Ext^1_H(N, \tau^{-1}B) \to
\Ext^1_H(N, \tau^{-1}M) \to 0.$$

\noindent {\bf Case (III):} \\
Finally, we consider the case where both $M$ and $M^{\ast}$ are
regular modules. By Lemma~\ref{lem3.3} (and arguing as in
Lemma~\ref{lem3.4})
we see that it is enough to show that either \\
(a) Any $H$-map $N\amalg \tau N\to M$ factors through $g$ in the sequence
\begin{equation} \label{modulesequence}
0\to M^{\ast}\overset{f}{\to} B\overset{g}{\to} M\to 0,
\end{equation}
or \\
(b) For any non-zero map $N\amalg \tau^{-1}N\to M$, the
composition with $h$ is non-zero, where $h$ is the map in the triangle:
$$\tau^{-1}M^{\ast} \to M\overset{h}{\to} B' \overset{k}{\to} M^{\ast}.$$

Note that statement (a) is the same as
statement~\ref{helpthm3.13}(i). In order to do this, it is clear
that we can assume that there is a non-zero $H$-map from $N\amalg
\tau N$ to $M$ and that there is a non-zero map (in $\C$) from
$N\amalg \tau^{-1}N$ to $M$. From the former, it follows that $N$
and $M$ must lie in the same tube.

We first of all show the following claim:

\noindent {\bf Claim:} If there is a non-zero H-map from $N\amalg
\tau^{-1}N$ to $M$ whose composition with $h$ is zero then any
non-zero $H$-map from $N\amalg \tau N$ to $M$ factors through $g:
B\to M.$

{\bf Proof of claim:} If $\gamma:N\to M$ is a non-zero $H$-map and
the composition $h\circ \gamma$ is zero, then $\gamma$ lifts to
$\tau^{-1}M^{\ast}$. It follows that $\tau^{-1}M^{\ast}$ lies in the
rectangle spanned by $N$ and $M$ in the tube (see
Figure~\ref{f:regreg}). Hence $M^{\ast}$ lies in the rectangle spanned by
$\tau M$ and $\tau N$ in the tube (the dotted rectangle in
Figure~\ref{f:regreg}). The middle term in the short exact
sequence~\eqref{modulesequence} is the direct sum of two
indecomposable modules, $B_1$ and $B_2$, which are at the other
two corners of the rectangle spanned by $M$ and $M^{\ast}$. We see from
the structure of the tube that $\gamma$ factors through at least
one of $B_1$ and $B_2$ (using the fact that $M^{\ast}\not\simeq \tau N$).
It follows that $\gamma$ factors through
the map $B\overset{g}{\to} M$. Since $\Hom_H(N,M)\simeq k$, we see
that every $H$-map from $N$ to $M$ factors through
$B\overset{g}{\to} M$.

Since a non-zero $H$-map from $\tau N$ to $M$ factors through any $H$-map
from $N$ to $M$ (note that such maps, if they exist, are unique up to a
scalar multiple), it follows that any non-zero $H$-map from $\tau N$ to
$M$ factors through $B\overset{g}\to M$
also. Thus any non-zero $H$-map from $N\amalg \tau N$
to $M$ factors through $B\overset{g}{\to} M$.

If $\gamma:\tau^{-1}N\to M$ is a non-zero $H$-map and its
composition with $h$ is zero, then, as above, $\gamma$ must factor
through $B\overset{g}\to M$. Then any non-zero $H$-map from $N$ to
$M$ must factor through $\gamma$ and therefore, by the above,
through $B\overset{g}\to M$. Again arguing as above, we see that
any non-zero $H$-map from $\tau N$ to $M$ factors through
$B\overset{g}\to M$ and therefore that any non-zero $H$-map from
$N\amalg \tau N$ to $M$ factors through $B\overset{g}\to M$. The
claim now follows.

Since $M$ is regular, Lemma~\ref{helpthm3.13} applies, and we see
that either~\ref{helpthm3.13}(i) or (ii) holds. If (i) holds, we
are done, as (a) above is the same as~\ref{helpthm3.13}(i).
If~\ref{helpthm3.13}(ii) holds, it follows that for any non-zero
$F$-map $N\amalg \tau^{-1}N \to M$, the composition $N\amalg
\tau^{-1}N \to M \to B'$ is nonzero.

By the claim above, we know that if there is a non-zero $H$-map
$N\amalg \tau^{-1}N \to M$ whose composition with $h$ is zero
then (a) holds, and we are done. So we are left with the case in
which for any non-zero $H$-map or $F$-map from $N\amalg
\tau^{-1}N$ to $M$ the composition with $h$ is non-zero.

We claim that this means that for any non-zero map from $N\amalg
\tau^{-1}N$ to $M$, the composition $N\amalg \tau^{-1}N \to M\to
B'$ is non-zero. We write maps in $\C$ between indecomposable
objects as pairs $(f,f')$ where $f$ is an $H$-map and $f'$ is an
$F$-map. Let $(h,h'):N\to M$ be a non-zero map such that the
composition with $M\to B'$ is zero. Then there are maps
$(t,t'):N\rightarrow \tau^{-1}M^{\ast}$ and $(s,s'):\tau^{-1}M^{\ast}
\rightarrow M$ such that $(h,h')=(s,s')\circ (t,t')$. Then $h=st$,
so that $h:N\to M$ composed with $M\to B'$ is zero; hence $h=0$
(as $h$ is a $H$-map). Then $h'$ composed with $M\to B'$ is zero,
so $h'=0$ (as $h'$ is an $F$-map). It follows that $(h,h')=0$,
a contradiction. See Figure~\ref{f:HFargument}.

We see that the map from $(N\amalg \tau^{-1}N,M)_{\C}$ to
$(N\amalg \tau^{-1}N,B')_{\C}$ is injective, and thus that (b)
above holds.

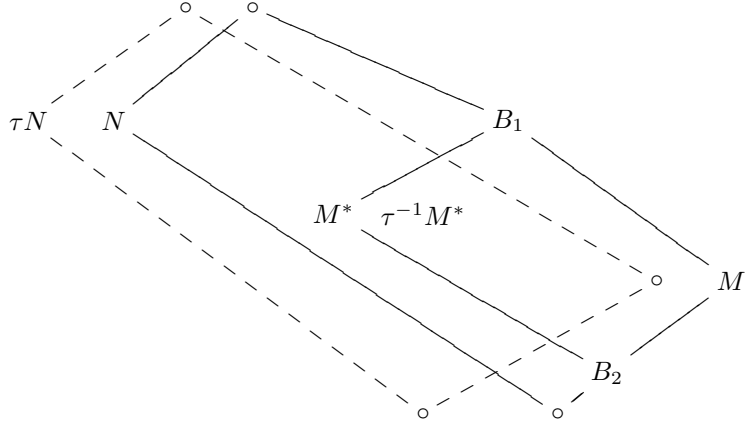
\begin{figure}
$$
\xy
\xymatrix@C=4pt@R=3pt{
&&&& \circ \ar@{--}[ddddllll] \ar@{--}[dddddddddrrrrrrrrr] &&
\circ \ar@{-}[ddddllll] \ar@{-}[ddddrrrr] &&&&&&&& \\
&&&&&&&&&&&&&&\\
&&&&&&&&&&&&&&\\
&&&&&&&&&&&&&& \\
\tau N \ar@{--}[dddddddddrrrrrrrrr] && N \ar@{-}[dddddddddrrrrrrrrr] &&&&&&&&B_1
\ar@{-}[dddll] \ar@{-}[dddddrrrrr] &&&& \\
&&&&&&&&&&&&&&&\\
&&&&&&&&&&&&&&\\
&&&&&&&& M^{\ast} \ar@{-}[dddddrrrr] & \tau^{-1}M^{\ast} &&&&&\\
&&&&&&&&&&&&&&\\
&&&&&&&&&&&&& \circ \ar@{--}[ddddllll] && M \ar@{-}[dddlll] \\
&&&&&&&&&&&&&& \\
&&&&&&&&&&&&&&& \\
&&&&&&&&&&&&B_2 \ar@{-}[dl] &&& \\
&&&&&&&&&\circ && \circ &&& \\
&&&&&&&&&&&&&&
}
\endxy
$$
\caption{Proof of claim in Theorem~\ref{thm3.13}}
\label{f:regreg}
\end{figure}

\begin{figure}
$$
\xy
\xymatrix{
& N \ar[dl]_{(t,t')} \ar[d]^{(h,h')} \\
\tau^{-1}M^{\ast} \ar[r]_{(s,s')} & M \ar[r] & B' 
}
\endxy
$$
\caption{End of proof of Theorem~\ref{thm3.13}}
\label{f:HFargument}
\end{figure}

We can now conclude that if $M \not \simeq \tau N \not \simeq M^{\ast}$, then 
\begin{multline*} \dim\Hom_{\C}(N,M)+\dim\Hom_{\C}(N,M^{\ast}) 
\\ = \max(\dim\Hom_{\C}(N,B),\dim\Hom_{\C}(N,B')). \end{multline*}
Thus $N$ is exchange compatible, and
we are done.
\end{proof}

\begin{cor}
Let $H$ be a tame hereditary algebra, and let $N$ be an indecomposable exceptional object in $\C$.
Then $N$ is exchange compatible if and only if $\End_{\C}(N)\simeq k$.
\end{cor}

\begin{proof}
This follows immediately from Theorem~\ref{thm3.13} and Corollary~\ref{thm3.12}
together with Proposition~\ref{p:endnotk}.
\end{proof}

\noindent {\bf Proof of Theorem~\ref{main2}:} By Theorem~\ref{t:all},
Theorem~\ref{main2}(a) is equivalent to the statement: \\
(a') All indecomposable direct summands $T_i$ of $T$ are exchange
compatible. \\
By Lemma~\ref{lem3.11}, \ref{main2}(c) implies \ref{main2}(b); by
Theorem~\ref{thm3.13}, \ref{main2}(b) implies (a'), and by
Corollary~\ref{c:endpropc}, (a') implies \ref{main2}(c). The proof is
complete.~$\Box$.

\section{Examples}
In this final section we give two examples to illustrate our main result.

\begin{example} \rm \label{e:A3}
We consider a Dynkin quiver $Q$ of type $A_3$, and choose a cluster
$\{y_1, y_2, y_3 \}$ with corresponding cyclic quiver $\Gamma$: 
$
\xy
\xymatrix@C=8pt@R=5pt{
1 \ar[rr] & & 2 \ar[ddl] \\
& & \\
& 3 \ar[uul] &
}
\endxy
$.

If we choose a cluster-tilting object $T$ where the quiver of $\End_{\C}(T)^{\op}$ is
$Q'$ and compute all cluster variables with respect to the initial cluster
$\{y_1, y_2, y_3 \}$, we get the picture shown in Figure~\ref{f:A3example},
where in the AR-quiver we have represented the indecomposable objects by
their corresponding cluster variables.

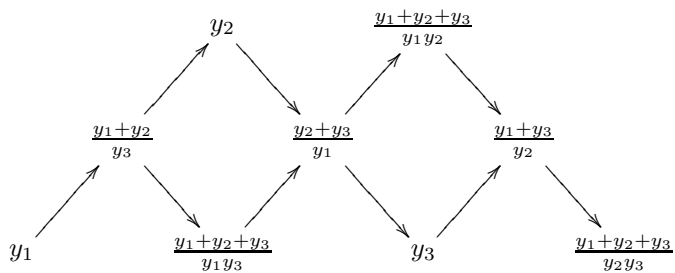
\begin{figure}
$$
\xy
\xymatrix@!@C=-7pt@R=-2pt{
& & y_2 \ar[dr] & & \frac{y_1 + y_2 + y_3}{y_1 y_2} \ar[dr] &  & \\  
& \frac{y_1 + y_2}{y_3} \ar[ur] \ar[dr] & & \frac{y_2 + y_3}{y_1} \ar[ur] \ar[dr] & & \frac{y_1 + y_3}{y_2} \ar[dr] & \\
y_1 \ar[ur] & & \frac{y_1 + y_2 + y_3}{y_1 y_3} \ar[ur] & &  y_3 \ar[ur] & &
\frac{y_1 + y_2 +y_3}{y_2 y_3}
}
\endxy
$$
\caption{Cluster variables in Example~\ref{e:A3}}
\label{f:A3example}
\end{figure}
It is easily verified here that the denominators are all of the
prescribed form, as promised by Theorem~\ref{main1}.
\end{example}

\begin{example} \rm \label{e:A2tilde}
In the next example we consider the quiver $Q$ shown in
Figure~\ref{f:A2tilde}(a). This is mutation equivalent to the quiver
$\Gamma$ shown in Figure~\ref{f:A2tilde}(b).
Then the AR-quiver of $kQ$ has a tube of rank two with quasisimple modules
given by the simple module $S_2$ and $M=\tau S_2$, a module with composition
factors $S_1$ and $S_3$.
We choose a cluster-tilting object $T$ such that the quiver of the cluster-tilted
algebra $\End_{C}(T)^{\op}$  is $\Gamma$, for example
$T = T_1 \amalg T_2 \amalg T_3 = P_1 \amalg M \amalg P_3$.
Then $\tau T = \tau P_1 \amalg S_2 \amalg  \tau P_3$.

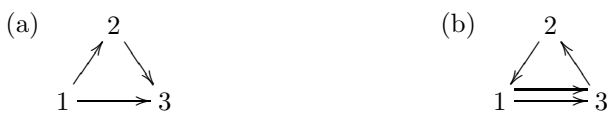
\begin{figure}
\begin{minipage}[b]{5cm} 
\centering
(a)
$
\xy
\xymatrix@C=8pt@R=5pt{
& 2 \ar[ddr] & \\
& & \\
1 \ar[uur] \ar[rr] & & 3
}
\endxy
$ 
\end{minipage}
\hspace{0.5cm} 
\begin{minipage}[b]{5cm}
\centering
(b)
$
\xy
\xymatrix@C=8pt@R=5pt{
& 2 \ar[ddl] & \\
& & \\
1 \ar@<1ex>[rr] \ar[rr] & & 3 \ar[uul] 
}
\endxy
$ 
\end{minipage}
\caption{The quivers $Q$ and $\Gamma$ in Example~\ref{e:A2tilde}.}
\label{f:A2tilde}
\end{figure}

Let $y_i$ be the cluster variable corresponding to $\tau T_i$ for $i=1,2,3$.
Then there is a seed $(\{y_1,y_2,y_3\},\Gamma)$ of $\mathcal{A}(Q)$.
The cluster variables corresponding to some of the objects
in the preprojective/preinjective component of the AR-quiver of
$\C_{kQ}$ are shown in Figure~\ref{f:A2tildevariables}.
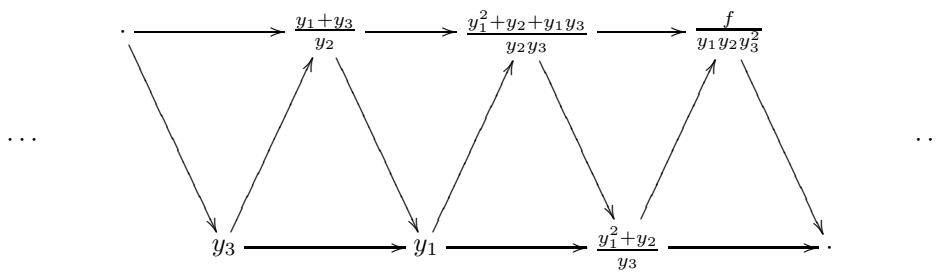
\begin{figure}
$$
\xy
\xymatrix@!@C=-15pt@R=-12pt{
& \cdot \ar[rr] \ar[ddr] & & \frac{y_1 + y_3}{y_2} \ar[rr] \ar[ddr]  & & \frac{y_1^2 + y_2 + y_1 y_3}{y_2 y_3} \ar[rr] \ar[ddr] & & 
\frac{f}{y_1 y_2 y_3^2}\ar[ddr] & &  \\
\cdots & & & & & & & & & \cdots \\
& & y_3 \ar[rr] \ar[uur] & & y_1 \ar[rr] \ar[uur] & & \frac{y_1^2 + y_2}{y_3}
\ar[rr] \ar[uur] & & \cdot & 
}
\endxy
$$
\caption{Some cluster variables in type $\widetilde{A_2}$}
\label{f:A2tildevariables}
\end{figure}

The cluster variables corresponding to the two quasisimple modules in the
tube of rank two are $y_2$, corresponding to
$\tau M=S_2$, and $\frac{(y_1+y_3)^2+y_2}{y_1y_2y_3}$ corresponding to $M$.
However $\Hom_{\C}(T_2,M) = \End_{\C}(M)$ is easily seen to be two
dimensional, while $y_2$ appears with multiplicity one in the denominator
of the cluster variable corresponding to $M$.

We remark that, using the above counter-example, C.-J. Fu and
B. Keller have discovered a counter-example
(see ~\cite{fukeller}) to Conjecture 7.17 of~\cite{fominzelevinsky07}.
See also \cite{ci}.
\end{example}


\begin{thebibliography}{99}

\bibitem[BFZ]{bfz05}
A. Berenstein, S. Fomin and A. Zelevinsky,
\emph{Cluster algebras III. Upper bounds and double Bruhat cells},
Duke Math. J. \textbf{126} (2005), no. 1, 1--52.

\bibitem[B]{bong}
K. Bongartz,
\emph{Tilted algebras}, 
Representations of algebras (Puebla, 1980), Lecture Notes in Math., 903, Springer, Berlin-New York, 1981, (26--38).

\bibitem[BCKMRT]{bckmrt}
A. B. Buan, P. Caldero, B. Keller, R. J. Marsh, I. Reiten and G. Todorov,
\emph{Appendix to Clusters and seeds in acyclic cluster algebras}, 
Proc. Amer. Math. Soc. 135, No. 10 (2007), 3049--3060.

\bibitem[BMR1]{bmrB}
A. B. Buan, R. Marsh and I. Reiten,
\emph{Cluster-tilted algebras},
Trans. Amer. Math. Soc., 359, no. 1 (2007), 323--332.

\bibitem[BMR2]{bmrA}
A. B. Buan, R. J. Marsh and I. Reiten,
\emph{Cluster mutation via quiver representations},
preprint arxiv:math.RT/0412077v2, 2004, to appear in Comment. Math. Helv.

\bibitem[BMRRT]{bmrrt06}
A. B. Buan, R. J. Marsh, M. Reineke, I. Reiten and G. Todorov,
\emph{Tilting theory and cluster combinatorics},
Advances in Mathematics \textbf{204} (2) (2006), 572--618.

\bibitem[BMRT]{bmrt}
A. B. Buan, R. J. Marsh, I. Reiten and G. Todorov,
\emph{Clusters and seeds in acyclic cluster algebras},
Proc. Amer. Math. Soc. 135, No. 10 (2007), 3049--3060.

\bibitem[CCS1]{ccs06}
P. Caldero, F. Chapoton and R. Schiffler
\emph{Quivers with relations arising from clusters ($A_n$ case)},
Transactions of the American Mathematical Society \textbf{358} (2006),
1347--1364.

\bibitem[CCS2]{ccs}
P. Caldero, F. Chapoton and R. Schiffler
\emph{Quivers with relations and cluster-tilted algebras},
Algebras and Representation Theory, 9, No. 4 , (2006), 359--376.

\bibitem[CK1]{calderokellerA}
P. Caldero and B. Keller,
\emph{From triangulated categories to cluster algebras},
preprint math.RT/0506018 (2005), to appear in Invent. Math.

\bibitem[CK2]{calderokellerB}
P. Caldero and B. Keller,
\emph{From triangulated categories to cluster algebras II},
Ann. Sci. Ecole Norm. Sup, 4eme serie, 39, (2006), 983--1009.

\bibitem[C]{ci} G. Cerulli Irelli, 
\emph{PhD-thesis in preparation: Structural theory of rank three cluster algebras of affine type}.

\bibitem[FZ1]{fominzelevinsky02}
S. Fomin and A. Zelevinsky,
\emph{Cluster algebras I: Foundations},
J. Amer. Math. Soc. \textbf{15} (2002), no. 2, 497--529.

\bibitem[FZ2]{fominzelevinsky03}
S. Fomin and A. Zelevinsky,
\emph{Cluster Algebras II: Finite type classification},
Invent. Math. \textbf{154}(1) (2003), 63--121.

\bibitem[FZ3]{fominzelevinsky07}
S. Fomin and A. Zelevinsky,
\emph{Cluster Algebras IV: Coefficients},
Compositio Mathematica \textbf{143} (2007), 112-164.

\bibitem[FK]{fukeller}
C. Fu and B. Keller,
\emph{On cluster algebras with coefficients and 2-Calabi-Yau categories},
preprint arxiv:0710.3152v1 [math.RT], 2007.

\bibitem[H]{happel}
D. Happel,
\emph{Triangulated categories in the representation theory of finite-dimensional algebras}, 
London Mathematical Society Lecture Note Series, 119. Cambridge University Press, Cambridge (1988).

\bibitem[HU]{hu89}
D. Happel and L. Unger,
\emph{Almost complete tilting modules},
Proc. Amer. Math. Soc. \textbf{107} (3) (1989), 603--610.

\bibitem[K1]{kern} O. Kerner, \emph{Stable components of wild tilted algebras}, J. Algebra 142 (1991), no. 1, 37--57.

\bibitem[K2]{kern2} O. Kerner, \emph{Representations of wild quivers}, 
Representation theory of algebras and related topics (Mexico City, 1994), 65--107, 
CMS Conf. Proc., 19, Amer. Math. Soc., Providence, RI, (1996).

\bibitem[MRZ]{mrz03}
R. Marsh, M. Reineke and A. Zelevinsky,
\emph{Generalized associahedra via quiver representations},
Trans. Amer. Math. Soc. 355 (2003), no. 1, 4171--4186.

\bibitem[R1]{rin}
C. M. Ringel,
\emph{Tame algebras and integral quadratic forms},
Springer Lecture Notes in Mathematics 1099 (1984).

\bibitem[Ri2]{rin2}
C. M. Ringel,
\emph{The regular components of the Auslander-Reiten quiver of a tilted algebra}, 
A Chinese summary appears in Chinese Ann. Math. Ser. A 9 (1988), no. 1, 102. Chinese Ann. Math. Ser. B 9 (1988), no. 1, 1--18.

\bibitem[RT]{reitentodorov}
I. Reiten and G. Todorov, unpublished.

\bibitem[S]{strauss} 
H. Strauss, \emph{On the perpendicular category of a partial tilting module} J. Algebra 144 (1991), no. 1, 43--66.

\end{thebibliography}
\end{document}